\documentclass[fleqn]{mat01}
\usepackage{times,mathtimy,amssymb,latexsym,rangecite}
\begin{document}

\setcounter{page}{233} \firstpage{233}

\font\xxx=msam10 at 10pt
\def\blacksquare{\mbox{\xxx{\char'245}\ \,}}

\newtheorem{theore}{Theorem}
\renewcommand\thetheore{\arabic{section}.\arabic{theore}}
\newtheorem{theor}[theore]{\bf Theorem}
\newtheorem{defn}[theore]{\rm DEFINITION}
\newtheorem{lem}[theore]{Lemma}
\newtheorem{remar}[theore]{Remark}
\newtheorem{pot}[theore]{Proof of Theorem}
\newtheorem{pol}[theore]{Proof of Lemma}
\newtheorem{case}{Case}
\newtheorem{step}{Step}

\renewcommand\theequation{\arabic{section}.\arabic{equation}}

\title{Infinite dimensional differential games with hybrid controls}

\markboth{A J Shaiju and Sheetal Dharmatti}{Infinite dimensional
differential games with hybrid controls}

\author{A J SHAIJU and SHEETAL DHARMATTI$^{*,\dagger}$}

\address{
Department of Mathematics, Indian Institute of Science,
Bangalore~560~012, India\\
\noindent $^*$TIFR Centre, IISc Campus, Bangalore~560~012, India\\
\noindent E-mail: shaiju@math.tifrbng.res.in; sheetal@math.iisc.ernet.in\\
\noindent $^{\dagger}$Corresponding author.}

\volume{117}

\mon{May}

\parts{2}

\pubyear{2007}

\Date{MS received 24 September 2005; revised 23 August 2006}

\begin{abstract}
A two-person zero-sum infinite dimensional differential game of
infinite duration with discounted payoff involving hybrid controls
is studied. The minimizing player is allowed to take continuous,
switching and impulse controls whereas the maximizing player is
allowed to take continuous and switching controls. By taking
strategies in the sense of Elliott--Kalton, we prove the existence
of value and characterize it as the unique viscosity solution of
the associated system of quasi-variational inequalities.
\end{abstract}


\maketitle

\section{Introduction and preliminaries}

The study of differential games with Elliott--Kalton strategies in
the viscosity solution framework is initiated by Evans and
Souganidis \cite{ES} where both players are allowed to take
continuous controls. Differential games where both players use
switching controls are studied by Yong \cite{Y1,Y2}. In \cite{Y3},
differential games involving impulse controls are considered; one
player is using continuous controls whereas the other uses impulse
control. In the final section of \cite{Y3}, the author mentions
that by using the ideas and techniques of the previous sections
one can study differential games where one player uses continuous,
switching and impulse controls and the other player uses
continuous and switching controls. The uniqueness result for the
associated system of quasi-variational inequalities (SQVI) with
bilateral constraints is said to hold under suitable non-zero loop
switching-cost condition and cheaper switching condition. In all
the above references, the state space is a finite-dimensional
Euclidean space.

The infinite dimensional analogue of \cite{ES} is studied by Kocan
{\it et~al} \cite{KSS}, where the authors prove the existence of
value and characterize the value function as the unique viscosity
solution (in the sense of \cite{CLVI}) of the associated
Hamilton--Jacobi--Isaacs equation.

In this paper, we study a two-person zero-sum differential game in
a Hilbert space where the minimizer (player~2) uses three types of
controls: continuous, switching and impulse. The maximizer
(player~1) uses continuous and switching controls. We first prove
dynamic programming principle (DPP) for this problem. Using DPP,
we prove that the lower and upper value functions are `approximate
solutions' of the associated SQVI in the viscosity sense
\cite{CLVI}. Finally we establish the existence of the value by
proving a uniqueness theorem for SQVI. We obtain our results
without any assumption like non-zero loop switching-cost condition
and/or cheaper switching-cost condition on the cost functions.
This will be further explained in the concluding section. Thus
this paper not only generalises the results of \cite{Y3} to the
infinite dimensional state space, it obtains the main result under
fairly general conditions as well.

The rest of the paper is organized as follows. We set up necessary
notations and assumptions in the remaining part of this section.
The statement of the main result is also given at the end of this
introductory section. The DPP is proved in \S2. In this section we
also show that the lower/upper value function is an `approximate
viscosity solution' of SQVI. Section~3 is devoted to the proof of
the main uniqueness result for SQVI and the existence of value. We
conclude the paper in \S4 with a few remarks.

We first describe the notations and basic assumptions. The state
space is a separable Hilbert space ${\mathbb{E}}$. The continuous
control set for player $i$, $i=1,2$, is $U^i$, a compact metric
space. The set $D^i=\{ d^i_1,\dots ,d^i_{m_i} \}$; $i=1,2$; is the
switching control set for player $i$. The impulse control set for
the player 2 is $K$, a closed and convex subset of the state space
${\mathbb{E}}$. The space of all $U^i$-valued measurable maps on
$[0,\infty)$ is the continuous control space for player $i$ and is
denoted by ${\mathcal{U}}^i$:
\begin{equation*}
{\mathcal{U}}^i = \{ u^i\hbox{\rm :}\ [0, \infty) \rightarrow U^i
| u^i ~\mbox{ measurable}\}.
\end{equation*}
By ${\mathcal{U}^i}[0, t]$ we mean the space of all $ U^i$-valued
measurable maps on $[0, t]$ that is,
\begin{equation*}
{\mathcal{U}}^i[0, t] = \{ u^i\hbox{\rm :}\ [0, t] \rightarrow U^i
| u^i ~\mbox{ measurable}\}.
\end{equation*}
The switching control space ${\mathcal{D}}^i$ for player $i$ and
the impulse control space ${\mathcal{K}}$ for player 2 are defined
as follows:
\begin{align*}
{\mathcal{D}}^i &= \Bigg\{ d^i(\cdot)=\sum_{j\geq 1}d^i_{j-1}
\chi_{[\theta^i_{j-1},\theta^i_j)}(\cdot) : d^i_j \in D^i,
(\theta^i_j)\subset [0,\infty],\\[.2pc]
&\quad \ \ \theta^i_0=0,
(\theta^i_j) \uparrow \infty , d^i_{j-1} \neq d^i_j \mbox{ if } \theta^i_j < \infty \Bigg\},\\[.3pc]
{\mathcal{K}} &= \Bigg\{ \xi(\cdot)=\sum_{j\geq 0} \xi_j
\chi_{[\tau_j,\infty]}(\cdot) : \xi_j \in K, (\tau_j)\subset
[0,\infty], (\tau_j) \uparrow \infty \Bigg\}.
\end{align*}

An impulse control $\xi(\cdot)=\sum_{j\geq 0} \xi_j
\chi_{[\tau_j,\infty]}(\cdot)$, consists of the impulse times
$\tau_j$'s and impulse vectors $\xi_j$'s. We use the notation
\begin{equation*}
(\xi)_{1,j} = \tau_j ~~\mbox{ and }~~ (\xi)_{2,j} = \xi_j .
\end{equation*}

Similarly for switching controls $d^1(\cdot)$ and $d^2(\cdot)$ we
write
\begin{align*}
(d^1)_{1, j} &= \theta^1_j\quad \mbox{ and }~~ (d^1)_{2, j} = d^1_{j},\\[.3pc]
(d^2)_{1, j} &= \theta^2_j\quad \mbox{ and }~~ (d^2)_{2, j} =
d^2_{j}.
\end{align*}

Now we describe the dynamics and cost functions involved in the
game. To this end, let ${\mathcal{C}}^1 = {\mathcal{U}}^1 \times
{\mathcal{D}}^1 $ and ${\mathcal{C}}^2 = {\mathcal{U}}^2 \times
{\mathcal{D}}^2 \times {\mathcal{K}}$. For
$(u^1(\cdot),d^1(\cdot)) \in {\mathcal{C}}^1 $ and
$(u^2(\cdot),d^2(\cdot),\xi(\cdot)) \in {\mathcal{C}}^2 $, the
corresponding state $y_x(\cdot)$ is governed by the following
controlled semilinear evolution equation in ${\mathbb{E}}$:
\begin{equation} \label{s1e1}
\dot{y}_x(t) + Ay_x(t) = f(y_x(t),u^1(t),d^1(t),u^2(t),d^2(t)) +
\dot{\xi}(t), \,\, y_x(0-)=x,
\end{equation}
where $f\hbox{\rm :}\ {\mathbb{E}} \times U^1 \times D^1 \times
U^2 \times D^2 \to {\mathbb{E}} $ and $-A$ is the generator of a
contraction semigroup $\{ S(t); t \geq 0 \}$ on
${\mathbb{E}}$.\vspace{.5pc}

\noindent {(A1)} We assume that the function $f$ is bounded,
continuous and for all $x,y \in {\mathbb{E}}$, $d^i \in D^i$, $u^i
\in U^i$,
\begin{equation}\label{s1e2}
\| f(x,u^1,d^1,u^2,d^2) - f(y,u^1,d^1,u^2,d^2) \| \leq L \| x-y
\|.
\end{equation}

Note that under the assumption (A1), for each $x \in
{\mathbb{E}}$, $ d^i(\cdot) \in {\mathcal{D}}^i $, $ u^i(\cdot)
\in {\mathcal{U}}^i$ and $ \xi(\cdot) \in {\mathcal{K}} $ there is
a unique mild solution $y_x(\cdot)$ of (\ref{s1e1}). This can be
concluded for example, from Corollary~2.11, chapter~4, page
number~109 of \cite{P}.

Let $k\hbox{\rm :}\ {\mathbb{E}} \times U^1 \times D^1 \times U^2
\times D^2 \to {\mathbb{R}} $ be the running cost function,
$c^i\hbox{\rm :}\ D^i \times D^i \to {\mathbb{R}}$ the switching
cost functions, and $l\hbox{\rm :}\ K \to {\mathbb{R}}$ the
impulse cost function.\vspace{.5pc}

\noindent {(A2)} We assume that the cost functions $k$, $c^i$, $l$
are nonnegative, bounded, continuous, and for all $x,y \in
{\mathbb{E}}$, $d^i \in D^i$, $u^i \in U^i$, $\xi_0 , \xi_1 \in
K$,
\begin{align}\label{s1e3}
\hskip -4pc \left.\begin{array}{@{}r@{\ }c@{\ }l} |
k(x,u^1,d^1,u^2,d^2) -
k(y,u^1,d^1,u^2,d^2) | &\leq& L \| x-y \|, \\[.7pc]
l(\xi_0 + \xi_1) &<& l(\xi_0) + l(\xi_1),
\, \forall \, \xi_0, \xi_1 \in K,\\[.7pc]
\displaystyle\lim_{|\xi | \rightarrow \infty } l(\xi) &=& \infty, \\[.9pc]
\displaystyle\inf_{d^i_1 \neq d^i_2} c^i(d^i_1,d^i_2) &=& c^i_0>0.
\end{array} \right\}\hskip -1pc\phantom{0}
\end{align}

\begin{remar}{\rm
The subadditivity condition $l(\xi_0 + \xi_1) <l(\xi_0) +
l(\xi_1)$ is needed to prove Lemma~\ref{s3i1} which, in turn, is
required to establish the uniqueness theorems (and hence the
existence of value for the game) in \S3. This condition makes sure
that, if an impulse $\xi_0$ is the best option at a particular
state $y_0$, then applying an impulse again is not a good option
for the new state $y_0+\xi_0$.}
\end{remar}

Let $\lambda >0$ be the discount parameter. The total discounted
cost functional $J_x\hbox{\rm :}\ {\mathcal{C}}^{1} \times
{\mathcal{C}}^{2} \to {\mathbb{R}}$ is given by
\begin{equation}\label{s1e4}
\hskip -4pc \left. \begin{array}{@{}r@{\ }c@{\ }l} J_x
[u^1(\cdot),d^1(\cdot),u^2(\cdot),d^2(\cdot),\xi(\cdot)]&
=&\displaystyle \int_{0}^{\infty} {\rm e}^{-\lambda t}
k(y_x(t),u^1(t),d^1(t),u^2(t),d^2(t)) ~{\rm
d}t\\[1.1pc]
&&\displaystyle - \sum_{j \geq 0} {\rm e}^{-\lambda \theta^1_j} c^1(d^1_{j-1},d^1_j) \\[1.3pc]
&&\displaystyle + \sum_{j \geq 0} {\rm e}^{-\lambda \theta^2_j}
c^2(d^2_{j-1},d^2_j) + \sum_{j \geq 1} {\rm e}^{-\lambda \tau_j}
l(\xi_j)
\end{array}\right\}.
\end{equation}

We next define the strategies for player 1 and player 2 in the
Elliott--Kalton framework. The strategy set $\Gamma$ for player 1
is the collection of all nonanticipating maps $\alpha$ from
${\mathcal{C}}^2$ to ${\mathcal{C}}^1$. The strategy set $\Delta$
for player 2 is the collection of all nonanticipating maps $\beta$
from ${\mathcal{C}}^1 $ to ${\mathcal{C}}^2$.

For a strategy $\beta$ of player 2 if
$\beta(u^1(\cdot),d^1(\cdot))=(u^2(\cdot),d^2(\cdot),\xi(\cdot))$,
then we write
\begin{align*}
\hskip -4pc \Pi_1 \beta (u^1(\cdot),d^1(\cdot))=u^2(\cdot), \;
\Pi_2 \beta (u^1(\cdot),d^1(\cdot))=d^2(\cdot)~ \mbox{ and }~
\Pi_3 \beta (u^1(\cdot),d^1(\cdot))=\xi(\cdot).
\end{align*}
That is, $\Pi_i$ is the projection on the $i$th component of the
map $\beta$. Similar notations are used for $\alpha
(u^2(\cdot),d^2(\cdot),\xi(\cdot))$ as well. Hence,
\begin{align*}
\Pi_1 \alpha (u^2(\cdot),d^2(\cdot), \xi(\cdot) )=u^1(\cdot) \quad
\mbox{and} \quad \Pi_2 \alpha (u^2(\cdot),d^2(\cdot), \xi(\cdot)
)=d^1(\cdot).
\end{align*}

Let ${\mathcal{D}}^{i,d^i}$ denote the set of all switching
controls for player $i$ starting at $d^i$. Then we define sets
\begin{equation*}
{\mathcal{C}}^{1,d^1} = {\mathcal{U}}^1 \times {\mathcal{D}}^{1,
d^1} \quad \mbox{and}\quad {\mathcal{C}}^{2,d^2}= {\mathcal{U}}^2
\times {\mathcal{D}}^{2, d^2} \times {\mathcal{K}}.
\end{equation*}
Let $\Delta^{d^2}$ denote the collection of all $\beta \in \Delta$
such that $\Pi_2 \beta (0)=d^2$ and $\Gamma^{d^1}$ be the
collection of all $ \alpha \in \Gamma$ such that $\Pi_2
\alpha(0)=d^1$.

Now using these strategies we define upper and lower value
functions associated with the game. Consider $J_x$ as defined in
(\ref{s1e4}). Let $J_x^{d^1,d^2}$ be the restriction of the cost
functional $J_x$ to ${\mathcal{C}}^{1,d^1} \times
{\mathcal{C}}^{2,d^2}. $ The upper and lower value functions are
defined respectively as follows:
\begin{align} \label{s1e5}
V_+^{d^1,d^2}(x) &= \sup_{\alpha \in \Gamma^{d^1}}
\inf_{{\mathcal{C}}^{2,d^2}}
J_x^{d^1,d^2}[\alpha(u^2(\cdot),d^2(\cdot),
\xi(\cdot)),u^2(\cdot),d^2(\cdot), \xi(\cdot)],\\[.3pc]\label{s1e6}
V_-^{d^1,d^2}(x) &= \inf_{\beta \in \Delta^{d^2}}
\sup_{{\mathcal{C}}^{1,d^1}}
J_x^{d^1,d^2}[u^1(\cdot),d^1(\cdot),\beta(u^1(\cdot),d^1(\cdot))].
\end{align}

Let $V_+ = \{ V_+^{d^1,d^2}\hbox{\rm :}\ (d^1,d^2) \in D^1 \times
D^2 \} $ and $V_- = \{ V_-^{d^1,d^2}\hbox{\rm :}\ (d^1,d^2) \in
D^1 \times D^2 \} $. If $V_+ \equiv V_- \equiv V$, then we say
that the differential game has a value and $V$ is referred to as
the value function.

Since all cost functions involved are bounded, value functions are
also bounded. In view of (A1) and (A2), the proof of uniform
continuity of $V_+$ and $V_-$ is routine. Hence both $V_+$ and
$V_-$ belong to ${\rm BUC}({\mathbb{E}};{\mathbb{R}}^{m_1 \times
m_2})$, the space of bounded uniformly continuous functions from
$\mathbb{E}$ to ${\mathbb{R}}^{m_1 \times m_2}$.

Now we describe the system of quasivariational inequalities (SQVI)
satisfied by upper and lower value functions and the definition of
viscosity solution in the sense of \cite{CLVI}.

For $x,p \in {\mathbb{E}}$, let
\begin{align}\label{s1e7}
\hskip -4pc H^{d^1,d^2}_-(x,p) &= \max_{u^2 \in U^2} \min_{u^1 \in
U^1} [
\langle -p,f(x,u^1,d^1,u^2,d^2) \rangle -k(x,u^1,d^1,u^2,d^2) ],\\[.3pc]\label{s1e8}
\hskip -4pc H^{d^1,d^2}_+(x,p) &= \min_{u^1 \in U^1} \max_{u^2 \in
U^2} [ \langle -p,f(x,u^1,d^1,u^2,d^2) \rangle
-k(x,u^1,d^1,u^2,d^2) ]
\end{align}
and for $V \in C({\mathbb{E}};{\mathbb{R}}^{m_1 \times m_2})$, let
\begin{align}\label{s1e9}
M_-^{d^1,d^2}[V](x) &= \min_{\bar{d}^2\neq d^2}
[V^{d^1,\bar{d}^2}(x)+ c^2(d^2,\bar{d}^2)], \\[.3pc]\label{s1e10}
M_+^{d^1,d^2}[V](x) &= \max_{\bar{d}^1\neq d^1}
[V^{\bar{d}^1,d^2}(x) - c^1(d^1,\bar{d}^1)],\\[.3pc]\label{s1e11}
N[V^{d^1,d^2}](x) &= \inf_{\xi \in K } [V^{d^1,{d}^2}(x+\xi)+
l(\xi)].
\end{align}

The HJI upper systems of equations associated with $V_-$ of the
hybrid differential game are as follows: for $(d^1,d^2) \in D^1
\times D^2,$
\begin{align*}
&\hskip -5.3pc \left. \begin{array} {@{}rcl} &&\min \, \{ \, \max
\, ( \lambda V^{d^1,d^2}+ \langle Ax, DV^{d^1,d^2} \rangle +
H^{d^1,d^2}_+(x,DV^{d^1,d^2}), \\[2mm]
&&\quad\, V^{d^1,d^2}-M_-^{d^1,d^2}[V], V^{d^1,d^2}-N[V^{d^1,d^2}]
\,), V^{d^1,d^2}- M^{d^1,d^2}_+[V] \, \} = 0
\end{array} \right\},\hskip -1pc\phantom{0} \tag{{\rm HJI1}+}\\[.6pc]
&\hskip -5.3pc \left. \begin{array} {@{}rcl} &&\max \, \{ \, \min
\, ( \lambda V^{d^1,d^2}+ \langle Ax, DV^{d^1,d^2} \rangle +
H^{d^1,d^2}_+(x,DV^{d^1,d^2}), \\[2mm]
&&\quad\, V^{d^1,d^2}-M_+^{d^1,d^2}[V] \,), \; \; V^{d^1,d^2}-
M^{d^1,d^2}_-[V], V^{d^1,d^2} - N[V^{d^1,d^2}] \, \} = 0
\end{array} \right\},\hskip -1pc\phantom{0} \tag{{\rm HJI2}+}
\end{align*}
where $ H^{d^1,d^2}_+, M^{d^1,d^2}_-, M^{d^1,d^2}_+ $ and $ N[
V^{d^1,d^2}] $ are as given by (\ref{s1e8}), (\ref{s1e9}),
(\ref{s1e10}) and (\ref{s1e11}) respectively.

Note here that for any real numbers $ a, b, c, d$, $ (a \vee b)
\wedge c \wedge d \leq (a \vee c \vee d ) \wedge b $.

If we replace $H^{d^1,d^2}_+$ in the above system of equations by
$H^{d^1,d^2}_-$, then we obtain the HJI lower system of equations
associated with $V_+$ and is denoted respectively by (HJI1$-$) and
(HJI2$-$).

If $V$ satisfies both (HJI1$+$) and (HJI2$+$), then we say that
$V$ satisfies (HJI$+$) and similarly if it satisfies both
(HJI1$-$) and (HJI2$-$), we say that $V$ satisfies (HJI$-$) .

Now let us recall the definition of viscosity solution (in the
sense of Crandall and Lions \cite{CLVI}). To this end, let
\begin{align*}
\hskip -4pc C^1({\mathbb{E}}) &= \{ \phi\hbox{\rm :}\ {\mathbb{E}}
\to {\mathbb{R}} \, |
\, \phi \mbox{ continuously differentiable}\},\\[.5pc]
\hskip -4pc {\rm Lip}({\mathbb{E}}) &= \{ \psi\hbox{\rm :}\
{\mathbb{E}} \to {\mathbb{R}}
\, | \, \psi \mbox{ Lipschitz continuous}\},\\[.5pc]
\hskip -4pc {\mathcal{T}} &= \{ \Phi \, | \, \Phi = \phi + \psi,
\, \phi \in
C^1({\mathbb{E}}), \, \psi \in {\rm Lip}({\mathbb{E}}) \},\\[.5pc]
\hskip -4pc D_A^+\Phi(x) &= \limsup_{\delta \downarrow 0, y \to x}
\frac{1}{\delta}[\Phi(y)-\Phi(S(\delta)y)],\\[.5pc]
\hskip -4pc D_A^-\Phi(x) &= \liminf_{\delta \downarrow 0, y \to x}
\frac{1}{\delta}[\Phi(y)-\Phi(S(\delta)y)],\\[.5pc]
\hskip -4pc H^{d^1,d^2}_{+,\bar{r}} (x,p) &= \sup_{\| q \| \leq r}
H^{d^1,d^2}_{+} (x,p+q); \quad H^{d^1,d^2}_{-,\bar{r}} (x,p) =
\sup_{\| q \| \leq r} H^{d^1,d^2}_{-} (x,p+q);\\[.5pc]
\hskip -4pc H^{d^1,d^2}_{+,\underline{r}} (x,p) &= \inf_{\| q \|
\leq r} H^{d^1,d^2}_{+} (x,p+q); \quad
H^{d^1,d^2}_{-,\underline{r}} (x,p) = \inf_{\| q \| \leq r}
H^{d^1,d^2}_{-} (x,p+q).
\end{align*}

\begin{defn}$\left.\right.$\vspace{.5pc}

\noindent{\rm A continuous function $V$ is a viscosity subsolution
of (HJI1$+$) if
\begin{align*}
\begin{array} {@{}rcl}
&&\hskip -5.1pc \min \, \{  \max \, ( \lambda
V^{d^1,d^2}(\hat{x})+ D_A^- \Phi(\hat{x}) +
H^{d^1,d^2}_{+,\underline{L(\psi)}}(\hat{x},D\phi(\hat{x})),
V^{d^1,d^2}(\hat{x})-M_-^{d^1,d^2}[V](\hat{x}), \\[.7pc]
&&\hskip -5.1pc \quad\,
V^{d^1,d^2}(\hat{x})-N[V^{d^1,d^2}](\hat{x})),
V^{d^1,d^2}(\hat{x})- M^{d^1,d^2}_+[V](\hat{x}) \} \leq 0,
\end{array}
\end{align*}
for any $\Phi \in {\mathcal{T}}$, $(d^1,d^2) \in D^1 \times D^2$
and local maximum $\hat{x}$ of $V^{d^1,d^2}-\Phi$.

A continuous function $V$ is a viscosity supersolution of
(HJI1$+$) if
\begin{align*}
\begin{array}{@{}rcl}
&&\hskip -5.1pc \min \, \{\max \, ( \lambda V^{d^1,d^2}(\hat{x})+
D_A^+ \Phi(\hat{x}) +
H^{d^1,d^2}_{+,\overline{L(\psi)}}(\hat{x},D\phi(\hat{x})),
V^{d^1,d^2}(\hat{x})-M_-^{d^1,d^2}[V](\hat{x}),\\[.7pc]
&&\hskip -5.1pc \quad\,
V^{d^1,d^2}(\hat{x})-N[V^{d^1,d^2}](\hat{x})),
V^{d^1,d^2}(\hat{x})- M^{d^1,d^2}_+[V](\hat{x})\} \geq 0,
\end{array}
\end{align*}
for any $\Phi \in {\mathcal{T}}$, $(d^1,d^2) \in D^1 \times D^2$
and local minimum $\hat{x}$ of $V^{d^1,d^2}-\Phi$.

If $V$ is both a subsolution and a supersolution of (HJI1$+$),
then we say that $V$ is a viscosity solution of (HJI1$+$).}
\end{defn}

\begin{defn}$\left.\right.$\vspace{.5pc}

\noindent{\rm A continuous function $V$ is an approximate
viscosity subsolution of (HJI1$+$) if for all $R>0$, there exists
a constant $C_R>0$ such that
\begin{align*}
\begin{array}{@{}rcl}
&&\hskip -5.1pc \min \, \{\max \, ( \lambda V^{d^1,d^2}(\hat{x})+
D_A^- \Phi(\hat{x}) + H^{d^1,d^2}_{+}(\hat{x},D\phi(\hat{x}) -C_R
L(\psi)),\\[.7pc]
&&\hskip -5.1pc \quad\,
V^{d^1,d^2}(\hat{x})-M_-^{d^1,d^2}[V](\hat{x}),
V^{d^1,d^2}(\hat{x})-N[V^{d^1,d^2}](\hat{x})),
V^{d^1,d^2}(\hat{x})- M^{d^1,d^2}_+[V](\hat{x}) \}\\[.7pc]
&&\hskip -5.1pc \quad\,\leq 0,
\end{array}
\end{align*}
for any $\Phi \in {\mathcal{T}}$, $(d^1,d^2) \in D^1 \times D^2$
and local maximum $\hat{x} \in B_R(0)$ of $V^{d^1,d^2}-\Phi$.

A continuous function $V$ is an approximate viscosity
supersolution of (HJI1$+$) if for all $R>0$, there exists a
constant $C_R>0$ such that
\begin{align*}
\begin{array}{@{}rcl}
&&\hskip -5.1pc \min \, \{ \max \, ( \lambda V^{d^1,d^2}(\hat{x})+
D_A^+ \Phi(\hat{x}) +
H^{d^1,d^2}_{+}(\hat{x},D\phi(\hat{x}) + C_R L(\psi)),\\[.7pc]
 &&\hskip -5.1pc \quad\, V^{d^1,d^2}(\hat{x})-M_-^{d^1,d^2}[V](\hat{x}),
V^{d^1,d^2}(\hat{x})-N[V^{d^1,d^2}](\hat{x})),
V^{d^1,d^2}(\hat{x})- M^{d^1,d^2}_+[V](\hat{x}) \}\\[.7pc]
&&\hskip -5.1pc \quad\, \geq 0,
\end{array}
\end{align*}
for any $\Phi \in {\mathcal{T}}$, $(d^1,d^2) \in D^1 \times D^2$
and local minimum $\hat{x} \in B_R(0)$ of $V^{d^1,d^2}-\Phi$.

If $V$ is both an approximate subsolution and an approximate
supersolution of (HJI1$+$), then we say that $V$ is an approximate
viscosity solution of (HJI1$+$).}
\end{defn}

In the above definitions, $L(\psi)$ is the Lipschitz constant of
$\psi$ and $B_R(0)$ is the closed ball of radius $R$ around the
origin.

\begin{remar}{\rm
One can easily prove that a viscosity solution is always an
approximate viscosity solution. For more details about the
approximate viscosity solution and its connections with other
notions of solutions, we refer to \cite{CLVI} and \cite{KSS}. In
the infinite dimensional set-up it is easier to establish that the
value functions are approximate viscosity solutions than to prove
that they are viscosity solutions. Therefore, as pointed out in
\cite{KSS}, the concept of approximate viscosity solution is used
as a vehicle to prove that the value functions are viscosity
solutions.}
\end{remar}

In the next section, we show that $V_-$ is an approximate
viscosity solution of (HJI$+$) and $V_+$ is an approximate
viscosity solution of (HJI$-$).

We say that the Isaacs min--max condition holds if
\begin{equation}\label{issacs}
H_-^{d^1,d^2} \equiv H_+^{d^1,d^2}\quad \mbox{ for all} \;
(d^1,d^2) \in D^1 \times D^2.
\end{equation}
Under this condition, the equations (HJI1$+$) and (HJI2$+$)
respectively coincide with (HJI1$-$) and (HJI2$-$). We now state
the main result of this paper; the proof will be worked out in
subsequent sections.

\begin{theor}[\!]\label{s1i1}
Assume {\rm (A1), (A2)} and the Isaacs min--max condition. Then
$V_- = V_+$ is the unique viscosity solution of {\rm (HJI}$+${\rm
)} {\rm (}or {\rm (HJI}$-${\rm ))} in ${\rm
BUC}({\mathbb{E}},{\mathbb{R}}^{m_1 \times m_2})$.
\end{theor}

\begin{remar}{\rm
The Isaacs min--max condition (\ref{issacs}) holds for the class
of problems where $f$ is of the form
\begin{equation*}
f(x,u^1,d^1,u^2,d^2) = f_1(x,u^1,d^1,d^2) + f_2(x,d^1,u^2,d^2).
\end{equation*}}
\end{remar}

\section{Dynamic programming principle}

\setcounter{equation}{0} \setcounter{theore}{0}

In this section, we first prove the dynamic programming principle
for the differential games with hybrid controls. We first state
the results. The proofs will be given later. Throughout this
section we assume (A1) and (A2).

\begin{lem}\label{s2i1}
For $(x, d^1, d^2) \in {\mathbb{E}} \times D^1 \times D^2$ and
$t>0,$
\begin{equation} \label{s2e1}
\left. \begin{array} {rcl} V_-^{d^1,d^2}(x) &= &\displaystyle
\inf_{\beta \in \Delta^{d^2} } \sup_{{\mathcal{C}}^{1,d^1}}
\Bigg[\int_0^t {\rm e}^{- \lambda s} k(y_x(s),u^1(s),d^1(s),\\[1.5pc]
&&\Pi_1 \beta (u^1,d^1)(s), \Pi_2 \beta (u^1, d^1) (s))~{\rm d}s\\[1pc]
&&\displaystyle - \sum_{\theta^1_j < t} {\rm e}^{-\lambda
\theta^1_j} c^1(d^1_{j-1},d^1_j)\\[1.5pc]
&&\displaystyle + \sum_{(\Pi_2 \beta)_{1,j} < t}
{\rm e}^{-\lambda (\Pi_2 \beta)_{1,j}} c^2((\Pi_2 \beta)_{2,j-1},(\Pi_2 \beta)_{2,j}) \\[1.5pc]
&&\displaystyle  + \sum_{(\Pi_3 \beta)_{1,j} < t}
{\rm e}^{-\lambda (\Pi_3 \beta)_{1,j}} l((\Pi_3 \beta)_{2,j}) \\[1.5pc]
&&\displaystyle  + {\rm e}^{-\lambda t} V_-^{d^1(t),\Pi_2 \beta
(t)}(y_x(t))\Bigg]
\end{array}\right\}.
\end{equation}
\end{lem}

\begin{lem}\label{s2i2}
For $(x, d^1, d^2) \in {\mathbb{E}} \times D^1 \times D^2 $ and
$t>0,$
\begin{align*}
\begin{array}{@{}rcl}
\hskip -4pc V_+^{d^1,d^2}(x) &= &\displaystyle\sup_{\alpha \in
\Gamma^{d^1} } \inf_{{\mathcal{C}}^{2,d^2}} \Bigg[\int_0^t {\rm
e}^{- \lambda s}
k(y_x(s),\alpha (u^2,d^2,\xi)(s),u^2(s),d^2(s),\xi(s))~{\rm
d}s\\[1.5pc]
\hskip -4pc &&\displaystyle -\sum_{(\Pi_2 \alpha)_{1,j} < t} {\rm
e}^{-\lambda (\Pi_2 \alpha)_{1,j}} c^1((\Pi_2 \alpha)_{1,j-1},
(\Pi_2 \alpha){1,j})\\[1.5pc]
\hskip -4pc &&\displaystyle + \sum_{\theta^2_j < t} {\rm e}^{-\lambda \theta^2_{j}} c^2(d^2_{j-1},d^2_{j})\\[1.5pc]
\hskip -4pc &&\displaystyle + \sum_{\tau_{j} < t} {\rm
e}^{-\lambda \tau_{j}} l(\xi_j) + {\rm e}^{-\lambda t} V_+^{\Pi_2
\alpha(t),d^2(t)}(y_x(t)) \Bigg].
\end{array}
\end{align*}
\end{lem}

\begin{lem}\label{s2i3}
The following results hold{\rm :}
\begin{enumerate}
\renewcommand\labelenumi{\rm (\roman{enumi})}
\leftskip .4pc
\item $M_+^{d^1,d^2}[V_-](x) \leq V_-^{d^1,d^2}(x) $.

\item $ V_-^{d^1,d^2}(x) \leq \min \{
M_-^{d^1,d^2}[V_-](x), N[V_-^{d^1,d^2}](x)\}$.

\item Let $(x,d^1,d^2)$ be such that strict inequality holds in {\rm (i)}.
Let $\bar{\beta} \in \Delta^{d^2}_0$. Then there exists $t_0 >0$
such that the following holds{\rm :}

For each $0 \leq t \leq t_0,$ there exists $u^{1,t}(\cdot) \in
{\mathcal{U}}^1[0,t]$ such that
\begin{align*}
\hskip -1.25pc V_-^{d^1,d^2}(x) - t^2 &\leq \int_0^t {\rm
e}^{-\lambda s}
k(y_x(s),u^{1,t}(s),d^1,\bar{\beta}(u^{1,t}(\cdot),d^1)(s))~{\rm
d}s\\[.4pc]
&\quad\, + {\rm e}^{-\lambda t}V_{-}^{d^1,d^2}(y_x(t)).
\end{align*}

\item Let $(x,d^1,d^2)$ be such that strict inequality holds in
{\rm (ii)}. Let $\bar{u}^1 \in U^1$. Then there exists $t_0 >0$
such that the following holds{\rm :}

For each $0 \!\leq\! t \!\leq\! t_0,$ there exists $\beta^t
\!\in\! \Delta^{d^2}$ with $(\Pi_2 \beta^t(\bar{u}^1,d^1))_{1,1},
(\Pi_3 \beta^t(\bar{u}^1,d^1))_{1,1} > t_0$ such that
\begin{align*}
\hskip -1.25pc V_-^{d^1,d^2}(x) + t^2 &\geq \int_0^t {\rm
e}^{-\lambda s} k(y_x(s),\bar{u}^1,d^1,
\beta^t(\bar{u}^1,d^1)(s))~{\rm d}s\\[.4pc]
&\quad\, + {\rm e}^{-\lambda t} V_-^{d^1,d^2}(y_x(t)) .
\end{align*}
\end{enumerate}
\end{lem}

\begin{lem}\label{s2i4}
The following results hold.

\begin{enumerate}
\renewcommand\labelenumi{\rm (\roman{enumi})}
\leftskip .4pc
\item $M_+^{d^1,d^2}[V_+](x) \leq V_+^{d^1,d^2}(x)$.

\item $ V_+^{d^1,d^2}(x) \leq \min \{ M_-^{d^1,d^2}[V_+](x),
N[V_+^{d^1,d^2}](x) \}$.

\item Let $(x,d^1,d^2)$ be such that strict inequality holds in
{\rm (ii)}. Let $\bar{\alpha} \in \Gamma^{d^1}_0$. Then there
exists $t_0 >0$ such that the following holds{\rm :}

For each $0 \leq t \leq t_0,$ there exists $u^{2,t}(\cdot) \in
{\mathcal{U}}^2[0,t]$ such that
\begin{align*}
\hskip -1.25pc V_+^{d^1,d^2}(x) + t^2 &\geq \int_0^t {\rm
e}^{-\lambda s} k(y_x(s),u^{2,t}(s),d^2,\bar{\alpha}(u^{2,t},d^2,
\xi^\infty)(s))~{\rm d}s\\[.4pc]
&\quad\, + {\rm e}^{-\lambda t}V_{+}^{d^1,d^2}(y_x(t)).
\end{align*}

\item Let $(x,d^1,d^2)$ be such that strict inequality
holds in {\rm (i)}. Let $\bar{u}^2 \in U^2$. Then there exists
$t_0 >0$ such that the following holds{\rm :}

For each $0 \leq t \leq t_0$, there exists $\alpha^t \in
\Gamma^{d^1}$ with $(\Pi_2
\alpha^t(\bar{u}^2,d^2,\xi^{\infty}))_{1,1} > t_0$ such that
\begin{align*}
\hskip -1.25pc V_+^{d^1,d^2}(x) - t^2 &\leq \int_0^t {\rm
e}^{-\lambda s} k(y_x(s),\bar{u}^2,d^2,\alpha^t(\bar{u}^2,d^2,
\xi^{\infty})(s))~{\rm d}s\\[.4pc]
&\quad\, + {\rm e}^{-\lambda t} V_+^{d^1,d^2}(y_x(t)) .
\end{align*}
\end{enumerate}
\end{lem}

We prove Lemmas~\ref{s2i1} and \ref{s2i3}. The proofs of
Lemmas~\ref{s2i2} and \ref{s2i4} are analogous.

\setcounter{theore}{0}
\begin{pol}{\rm
Let $ (x, d^1, d^2) \in {\mathbb{E}} \times D^1 \times D^2 $ and
$t > 0$. Let us denote the RHS of (\ref{s2e1}) by $W(x)$. Fix
$\epsilon >0$.

Let $\bar{\beta} \in \Delta^{d^2}$ be such that
\begin{align*}
\begin{array}{@{}rcl}
\hskip -4pc W(x) &\geq &\displaystyle\sup_{{\mathcal{C}}^{1,d^1}}
\Bigg[\int_0^t {\rm e}^{- \lambda s}
k(y_x(s),u^1(s),d^1(s),\Pi_1\bar{\beta}(u^1,d^1)(s),
\Pi_2\bar{\beta} (u^1, d^1) (s) )~{\rm d}s \\[1.2pc]
\hskip -4pc &&\displaystyle - \sum_{\theta^1_j < t} {\rm
e}^{-\lambda \theta^1_j} c^1(d^1_{j-1},d^1_j) +
\sum_{(\Pi_2\bar{\beta})_{1,j} < t} {\rm e}^{-\lambda
(\Pi_2\bar{\beta})_{1,j}}
c^2((\Pi_2\bar{\beta})_{2,j-1},(\Pi_2\bar{\beta})_{2,j}) \\[.7pc]
\hskip -4pc &&\displaystyle + \sum_{(\Pi_3\bar{\beta})_{1,j} < t}
{\rm e}^{-\lambda (\Pi_3\bar{\beta})_{1,j}}
l((\Pi_3\bar{\beta})_{2,j}) + {\rm e}^{-\lambda t}
V_-^{d^1(t),(\Pi_2\bar{\beta})(t)}(y_x(t))\Bigg] - \epsilon .
\end{array}
\end{align*}

By the definition of $V_-$, for each $(u^1(\cdot),d^1(\cdot)) \in
{\mathcal{C}}^{1,d^1}$, there exists
$\beta_{u^1(\cdot),d^1(\cdot)} \in \Delta^{(\Pi_2\bar{\beta})(t)}$
such that
\begin{align*}
&V_-^{d^1(t),(\Pi_2\bar{\beta})(t)}(y_x(t))\\[.3pc]
&\quad\, \geq J^{d^1(t),
((\Pi_2\bar{\beta})(t)}_{y_x(t)}[u^1(\cdot),d^1(\cdot),
\beta_{u^1(\cdot),d^1(\cdot)}(u^1(\cdot),d^1(\cdot))]-\epsilon.
\end{align*}

Define $\delta \in \Delta^{d^2}$ by
\begin{align*}
\delta(u^1(\cdot),d^1(\cdot))(s)= \left \{ \begin{array}{lcr}
\bar{\beta}(u^1(\cdot),d^1(\cdot))(s); & & s \leq t \\[.4pc]
\beta_{u^1(\cdot),d^1(\cdot)}(u^1(\cdot+t),d^1(\cdot+t))(s-t); & &
s>t
\end{array}\right..
\end{align*}

By change of variables, we get
\begin{align*}
\begin{array}{@{\hskip -5.1pc}rcl}
&&\displaystyle
J^{d^1(t),(\Pi_2\bar{\beta})(t)}_{y_x(t)}[u^1(\cdot +t),d^1(\cdot
+t), \beta_{u^1(\cdot),d^1(\cdot)}(u^1(\cdot +t),d^1(\cdot +t))] \\[1.3pc]
&&\displaystyle\quad\, =\int_t^{\infty} {\rm e}^{- \lambda \tau }
k(y_x(\tau),u^1(\tau),d^1(\tau), (\Pi_1\delta)(u^1,d^1)(\tau),
(\Pi_2\delta)( u^1, d^1) (\tau))~{\rm d}\tau \\[1.3pc]
&&\displaystyle\qquad\ -\sum_{\theta^1_j > t}
{\rm e}^{-\lambda \theta^1_j} c^1(d^1_{j-1},d^1_j) \\[1.5pc]
&&\displaystyle\qquad\ +\sum_{(\Pi_2\delta)_{1,j} > t} {\rm
e}^{-\lambda (\Pi_2\delta)_{1,j}}
c^2((\Pi_2\delta)_{2,j-1},(\Pi_2\delta)_{2,j}) \\[1.5pc]
&&\displaystyle\qquad\ +\sum_{(\Pi_3\delta)_{1,j} > t} {\rm
e}^{-\lambda (\Pi_3\delta)_{1,j}} l((\Pi_3\delta)_{2,j}).
\end{array}
\end{align*}

Substituting above in the inequality of
$V_-^{d^1(t),(\Pi_2\bar{\beta})(t)}(y_x(t))$ and then in the
inequality for $ W(x)$ will imply
\begin{equation*}
W(x) \geq J^{d^1,d^2}_x
[u^1(\cdot),d^1(\cdot),\delta(u^1,d^1)(\cdot)] - 2\epsilon.
\end{equation*}
This holds for all $(u^1(\cdot),d^1(\cdot)) \in
{\mathcal{C}}^{1,d^1}$ and hence $W(x) \geq V^{d^1,d^2}_-(x) - 2
\epsilon $. Since $\epsilon >0$ is arbitrary, we get $W(x) \geq
V^{d^1,d^2}_-(x)$.

We now prove the other type of inequality. Fix $\beta \in
\Delta^{d^2}$ and $\epsilon >0$. Choose
$(\bar{u}^1(\cdot),\bar{d}^1(\cdot)) \in {\mathcal{C}}^{1,d^1}$
such that
\begin{equation} \label{s2e2}
\left. \begin{array} {rcl} W(x) &\leq&\displaystyle\int_0^t {\rm
e}^{- \lambda s} k(y_x(s),\bar{u}^1(s),\bar{d}^1(s),
\beta(\bar{u}^1,\bar{d}^1)(s))~{\rm d}s \\[1.3pc]
&&\displaystyle- \sum_{\theta^1_j < t}
{\rm e}^{-\lambda \theta^1_j} c^1(d^1_{j-1},d^1_j) \\[1.6pc]
&&\displaystyle+ \sum_{(\Pi_2{\beta})_{1,j} < t} {\rm e}^{-\lambda
(\Pi_2{\beta})_{1,j}}
c^2((\Pi_2{\beta})_{2,j-1},(\Pi_2{\beta})_{2,j}) \\[1.6pc]
&&\displaystyle+ \sum_{(\Pi_3{\beta})_{1,j} < t}
{\rm e}^{-\lambda (\Pi_3{\beta})_{1,j}} l((\Pi_3{\beta})_{2,j}) \\[1.6pc]
&&\displaystyle+ {\rm e}^{-\lambda t} V_-^{\bar{d}^1(t),(\Pi_2
\beta)(t)}(y_x(t)) + \epsilon
\end{array} \right\}.
\end{equation}

Now for each $u^1(\cdot)$, define $\tilde{u}^1(\cdot)$ by
\begin{equation*}
\tilde{u}^1(s)= \left \{
\begin{array}{ll}
\bar{u}^1(s); & s\leq t \\[.5pc]
u^1(s-t); & s>t
\end{array} \right..
\end{equation*}

Similarly, for each $d^1(\cdot)$, we define $\tilde{d}^1(\cdot)$.
Let
\begin{equation*}
\hat{\beta}(u^1(\cdot),d^1(\cdot))(s)=\beta(\tilde{u}^1(\cdot),
\tilde{d}^1(\cdot))(s+t).
\end{equation*}

By the definition of $V_-$, we can choose $(u^1(\cdot),d^1(\cdot))
\in {\mathcal{C}}^{1,\bar{d}^1(t)}$ such that
\begin{align}\label{s2e3}
&V_-^{\bar{d}^1(t),(\Pi_2 \beta)(t)}(y_x(t))\nonumber\\[.4pc]
&\quad\, \leq J^{\bar{d}^1(t),(\Pi_2 \beta)(t)}_{y_x(t)}
[u^1(\cdot + t ),d^1(\cdot +t),\hat{\beta}(u^1,d^1)(\cdot +t)] +
\epsilon {\rm e}^{\lambda t}.
\end{align}

Now, combining (\ref{s2e2}) and (\ref{s2e3}), we get
\begin{align*}
\begin{array}{@{\hskip -4.1pc}rcl}
W(x) &\leq &\displaystyle \int_0^t {\rm e}^{- \lambda s}
k(y_x(s),\bar{u}^1(s),\bar{d}^1(s),\beta(\bar{u}^1,\bar{d}^1)(s))~{\rm
d}s - \sum_{\theta^1_j < t}
{\rm e}^{-\lambda \theta^1_j} c^1(d^1_{j-1},d^1_j) \\[1.5pc]
&&\displaystyle + \sum_{(\Pi_2{\beta})_{1,j} < t} {\rm
e}^{-\lambda \Pi_2{\beta}_{1,j}}
c^2((\Pi_2{\beta})_{2,j-1},(\Pi_2{\beta})_{2,j})\\[1.5pc]
&&\displaystyle + \sum_{\Pi_3{\beta}_{1,j} < t}
{\rm e}^{-\lambda (\Pi_3{\beta})_{1,j}} l((\Pi_3{\beta})_{2,j}) \\[1.5pc]
&&\displaystyle + {\rm e}^{- \lambda t} J^{\bar{d}^1(t),(\Pi_2
\beta)(t)}_{y_x(t)}
[u^1(\cdot),d^1(\cdot),\hat{\beta}(u^1(\cdot),d^1(\cdot))] +
2\epsilon.
\end{array}
\end{align*}

By change of variables, it follows that
\begin{equation*}
W(x) \leq J^{d^1,d^2}_x[\tilde{u}^1(\cdot),\tilde{d}^1(\cdot),
{\beta}(\tilde{u}^1,\tilde{d}^1)(\cdot)] + 2 \epsilon.
\end{equation*}

This holds for any $\beta \in \Delta^{d^2}$
and hence
\begin{equation*}
W(x) \leq V_-^{d^1,d^2}(x)+ 2 \epsilon .
\end{equation*}

The proof is now complete, since $\epsilon $ is arbitrary.}\hfill
$\blacksquare$
\end{pol}

\setcounter{theore}{2}
\begin{pol}{\rm
We first prove (i) and (ii). By the definition of $V^-$, for any
$\bar{d}^2 \neq d^2$,
\begin{equation*}
V_-^{d^1,d^2}(x) \leq V_-^{d^1,\bar{d}^2}(x) + c^2(d^2,\bar{d}^2).
\end{equation*}

From this we get
\begin{equation*}
V_-^{d^1,d^2}(x) \leq M^{d^1,d^2}_-[V_-](x).
\end{equation*}

The inequality
\begin{equation*}
V_-^{d^1,d^2}(x) \geq M^{d^1,d^2}_+[V_-](x)
\end{equation*}
can be proved in a similar fashion.

Clearly for any $\xi \in K$ and $\beta \in \Delta^{d^2} $,
\begin{equation*}
V_-^{d^1,d^2}(x) \leq \sup_{{\mathcal{C}}^{1,d^1}}
J_{x+\xi}^{d^1,d^2}[u^1 (\cdot),d^1 (\cdot), \beta(u^1,d^1)
(\cdot)] + l(\xi).
\end{equation*}

Taking infimum over $\beta \in \Delta^{d^2} $ and then over $\xi
\in K$, we obtain
\begin{equation*}
V_-^{d^1,d^2}(x) \leq N[V_{-}^{d^1,d^2}](x).
\end{equation*}

We now turn to the proof of (iii). By Lemma~\ref{s2i1}, for each
$t \geq 0$, there exists $(u^{1,t}(\cdot),d^{1,t}(\cdot)) \in
{\mathcal{C}}^{1,d^1}$ such that
\begin{align*}
\begin{array}{@{}rcl}
\hskip -4pc V_-^{d^1,d^2}(x) - t^2 & \leq &\displaystyle \int_0^t
{\rm e}^{-\lambda s} k(y_x(s),u^{1,t}(s),d^{1,t}(s),
\bar{\beta}(u^{1,t},d^{1,t})(s))~{\rm d}s \\[1.5pc]
\hskip -4pc &&\displaystyle - \sum_{\theta^{1,t}_j < t} {\rm
e}^{-\lambda \theta^{1,t}_j} c^1(d^{1,t}_{j-1},d^{1,t}_j ) + {\rm
e}^{-\lambda t} V_{-}^{d^{1,t}(t),(\Pi_2\bar{\beta})(t)}(y_x(t)).
\end{array}
\end{align*}

It is enough to show that, for some $t_0>0$, $\theta^{1,t}_1 \geq
t $ for all $0 \leq t \leq t_0$. If this does not happen, then
there would exist a sequence $t_n \downarrow 0$ such that
$\theta^{1,t_n}_1 < t_n$ for all $n$. This would imply that
\begin{align*}
\begin{array}{@{}rcl}
\hskip -4pc V_-^{d^1,d^2}(x) - (\theta^{1,t_n}_1)^2 & \leq
&\displaystyle \int_0^{\theta_1^{1,t_n}} {\rm e}^{-\lambda s}
k(y_x(s),u^{1,t_n}(s),d^{1,t_n}(s),
\bar{\beta}(u^{1,t_n},d^{1,t_n})(s))~{\rm d}s \\[1pc]
\hskip -4pc &&\displaystyle - {\rm e}^{-\lambda \theta^{1,t_n}_1}
c^1(d^1,d^{1,t_n}_1 ) + {\rm }{\rm e}^{-\lambda \theta_1^{1,t_n}}
V_{-}^{d^{1,t_n}_1,d^2}(y_x(\theta_1^{1,t_n})).
\end{array}
\end{align*}

We may assume that for all $n$, $d^{1,t_n}_1 =\bar{d}^1\neq d^1 $.
Now by letting $n \to \infty$ in the above inequality, we get
\begin{equation*}
\begin{array} {rcl}
V_-^{d^1,d^2}(x) & \leq & - c^1(d^1,\bar{d}^1 ) +
 V_{-}^{\bar{d}^1,d^2}(x) \\[.5pc]
&\leq& M_+^{d^1,d^2} [V_-](x).
\end{array}
\end{equation*}

This contradicts the hypothesis that strict inequality holds in
(i) and the proof of (iii) is now complete.

We next prove (iv). By Lemma~\ref{s2i1}, for each $t >0$, there
exists $\beta^t \in \Delta^{d^2} $ such that
\begin{align*}
\begin{array}{@{\hskip -4pc}rcl}
V_-^{d^1,d^2}(x) + t^2 & \geq &\displaystyle \int_0^t {\rm
e}^{-\lambda s} k(y_x(s),\bar{u}^1 ,d^1,
\beta^t(\bar{u}^1,d^1)(s))~{\rm d}s \!+
{\rm e}^{-\lambda t} V_-^{d^1,(\Pi_2\beta^t)(t)}(y_x(t)) \\[1.2pc]
&&\displaystyle + \sum_{(\Pi_2\beta^t)_{1,j} < t}
{\rm e}^{- \lambda \theta_j^{2,t}} c^2(d_{j-1}^{2,t},d_{j}^{2,t})\\[1.5pc]
&&\displaystyle + \sum_{(\Pi_3 \beta^t)_{1,j} < t} {\rm e}^{-
\lambda \tau_j^{t}} l(\xi_j^t).
\end{array}
\end{align*}

It is enough to show that, for some $t_0 >0$, $\min
(\theta_1^{2,t}, \tau_1^t)\geq t$ for all $0 \leq t \leq t_0$. If
this were not true, then (without any loss of generality) there
would be a sequence $t_n \downarrow 0 $ and two cases to consider.
In the first case, $\theta_1^{2,t_n} \leq \min (t_n,\tau_1^{t_n})$
whereas in the second case $\tau_1^{t_n} \leq \min
(t_n,\theta_1^{2,t_n})$. By dropping to a subsequence if necessary
and proceeding as in the proof of (iii), we get $V^{d^1,d^2}_-(x)
\geq M^{d^1,d^2}_-[V_-](x)$ in case 1 and $V^{d^1,d^2}_-(x) \geq
N[V_-^{d^1,d^2}](x)$ in case 2 respectively. This contradicts our
hypothesis that strict inequality holds in (ii) and the proof is
complete.}\hfill $\blacksquare$
\end{pol}

\setcounter{theore}{4}
\begin{remar}{\rm
From the proofs it is clear that, instead of the term $t^2$ in the
statement of Lemma~\ref{s2i3} and Lemma~\ref{s2i4}, we can take
any modulus $\rho(t)$.}
\end{remar}

We now state the following result (Corollary~4.9, \cite{KSS})
which is useful in proving that $V_-$ (resp. $V_+$) is a viscosity
solution of (HJI$+$) (resp. (HJI$-$)).

\begin{lem} \label{c4.9}
Let $\Phi \in {\mathcal{T}}$, $\hat{x} \in {\mathbb{E}}$ and
$D_A^- \Phi (\hat{x}) < \infty$. If $v(\cdot) \in
L^1([0,T];{\mathbb{E}})$ and $y(\cdot)$ solves
\begin{equation*}
\begin{array}{l}
\dot{y}(t)+ Ay (t)=v(t) ; \,\, \,\, 0\leq t \leq T, \\[.6pc]
y(0)=\hat{x},
\end{array}
\end{equation*}
then as $t \to 0$,
\begin{align*}
\begin{array} {@{}rcl@{}}
{\rm e}^{-\lambda t} \Phi(y(t)) \!-\! \Phi(\hat{x}) & \leq
&\displaystyle \frac{{\rm e}^{-\lambda t}-1}{\lambda}D_A^- \Phi
(\hat{x}) +
L(\psi) \int_0^t {\rm e}^{-\lambda s} \| v(s) \|~{\rm d}s\\[1.5pc]
&&\displaystyle + \int_0^t {\rm e}^{-\lambda s} [\langle
D\phi(y(s)), v(s) \rangle \!-\! \lambda \Phi(y(s))]~{\rm d}s \!+\!
\small{o}(t),
\end{array}
\end{align*}
uniformly for all $v(\cdot)$ uniformly integrable on $(0,T)$.
\end{lem}

We are now ready to prove that $V_-$ (resp. $V_+$) is an
approximate viscosity solution of (HJI$+$) (resp. (HJI$-$)).

\begin{theor}[\!] \label{s2i6}
The lower value function $V_-$ is an approximate viscosity
solution of {\rm (HJI}$+${\rm )} and the upper value function
$V_+$ is an approximate viscosity solution of {\rm (HJI}$-${\rm
)}.
\end{theor}

\begin{proof}
We prove that $V_-$ is an approximate viscosity solution of
(HJI$+$). The other part can be proved in an analogous manner.

Let $C_R$ be such that $\|
f(y_{x}(t),u^1(t),d^1(t),u^2(t),d^2(t),\xi(t)) \| \leq C_R$ for
all $0\leq t \leq 1$, $\| x\| \leq R$, $(u^1(\cdot),d^1(\cdot))
\in {\mathcal{C}}^1$ and $(u^2(\cdot),d^2(\cdot),\xi(\cdot)) \in
{\mathcal{C}}^2$.

We first prove that $V_-$ is an approximate subsolution of
(HJI1$+$). Let $(\hat{x},d^1,d^2) \in {\mathbb{E}} \times D^1
\times D^2 $ and $\Phi=\phi + \psi \in {\mathcal{T}}$ be such that
$V_-^{d^1,d^2}-\Phi$ has a local maximum at $\hat{x} \in B_R(0)$.
Without any loss of generality, we may assume that
$V_-^{d^1,d^2}(\hat{x})=\Phi(\hat{x})$. If
$V_-^{d^1,d^2}(\hat{x})=$ $M^{d^1,d^2}_+[V_-](\hat{x})$, then we
are done. Assume that $V_-^{d^1,d^2}(\hat{x}) >
M^{d^1,d^2}_+[V_-](\hat{x})$. It suffices to show that
\begin{equation*}
\lambda \Phi(\hat{x}) + D_A^-\Phi(\hat{x}) +
H^{d^1,d^2}_+(\hat{x},D\phi(\hat{x})) -C_R L(\psi)=: r \leq 0.
\end{equation*}

If possible, let $r>0$. This implies that for every $u^1 \in U^1$,
there exists $u^2=u^2(u^1)$ such that
\begin{align*}
&\lambda \Phi(\hat{x}) + D_A^-\Phi(\hat{x}) - \langle
D\phi(\hat{x}), f(\hat{x},u^1,d^1,u^2,d^2) \rangle\\[.4pc]
&\quad\, -k(\hat{x},u^1,d^1,u^2,d^2) -C_R L(\psi) \geq
\frac{r}{2}.
\end{align*}

By Proposition~3.2 of \cite{KSS}, for $t$ small enough, there
exists $\beta^t \in \Delta^{d^2}_0$ such that
\begin{align*}
&\lambda \Phi(\hat{x}) + D_A^-\Phi(\hat{x}) - \langle
D\phi(\hat{x}),f(\hat{x},u^1(s),d^1, \beta^t(u^1(\cdot),d^1)(s))\rangle\\[.4pc]
&\quad\, -k(\hat{x},u^1(s),d^1,\beta^t(u^1(\cdot),d^1)(s)) -C_R
L(\psi) \geq \frac{r}{2}
\end{align*}
for all $u^1(\cdot) \in {\mathcal{U}}^1[0,t]$ and a.e. $s
\in[0,t]$.

This yields
\begin{align*}
&\lambda \Phi(y_{\hat{x}}(s)) + D_A^-\Phi(\hat{x}) - \langle
D\phi(y_{\hat{x}}(s)),f(y_{\hat{x}}(s),u^1(s),d^1,\beta^t(u^1(\cdot),d^1)(s))
\rangle\\[.4pc]
&\quad\, -k(y_{\hat{x}}(s),u^1(s),d^1,\beta^t(u^1(\cdot),d^1)(s))
-C_R L(\psi) \geq \frac{r}{2}
\end{align*}
for a.e. $s \in [0,t]$ and $t$ small enough.

Multiplying throughout by ${\rm e}^{-\lambda s}$ and integrating
from $s=0$ to $s=t$, we get
\begin{equation} \label{1}
\left.\begin{array}{@{}rcl} &&\displaystyle\left[
D_A^-\Phi(\hat{x}) - C_R L(\psi)-\frac{r}{2}\right]
\left[\frac{{\rm e}^{-\lambda t}-1}{-\lambda} \right]\\[1.3pc]
&&\displaystyle \quad\, - \int_0^t {\rm e}^{-\lambda s}
k(y_{\hat{x}}(s),u^1(s),d^1,\beta^t(u^1(\cdot),d^1)(s))~{\rm
d}s\\[1.5pc]
&&\quad\, \displaystyle + \int_0^t {\rm e}^{-\lambda s}[\lambda
\Phi(y_{\hat{x}}(s))\\[1.5pc]
&&\quad\, - \langle D\phi(y_{\hat{x}}(s)),
f(y_{\hat{x}}(s),u^1(s),d^1,\beta^t(u^1(\cdot),d^1)(s)) \rangle ]
~{\rm d}s \geq 0
\end{array}\right\}.
\end{equation}

By Lemma~\ref{s2i3}(iii), for $t$ small enough, there exists
$u^{1,t}(\cdot) \in {\mathcal{U}}^1[0,t] $ such that
\begin{align} \label{2}
V_-^{d^1,d^2}(\hat{x}) - \small{o}(t) &\leq \int_0^t {\rm
e}^{-\lambda s}
k(y_{\hat{x}}(s),u^{1,t}(s),d^1,{\beta}^t(u^{1,t}(\cdot),d^1)(s))~{\rm
d}s\nonumber\\[.7pc]
&\quad\, + {\rm e}^{-\lambda t}V_{-}^{d^1,d^2}(y_{\hat{x}}(t)).
\end{align}

We now claim that $D_A^-\Phi(\hat{x}) < \infty$. We may take
$\Phi$ to be Lipschitz.
\begin{equation*}
\begin{array}{r@{\ }c@{\ }l}
\| \Phi(\hat{x})-\Phi(y_{\hat{x}}(\delta)) \| &\leq& \|
\Phi(\hat{x})-{\rm e}^{-\lambda \delta} \Phi(y_{\hat{x}}(\delta))
\| + ({\rm e}^{-\lambda \delta}-1) \| \Phi(y_{\hat{x}}(\delta)) \|\\[.7pc]
&\leq& \small{o}(\delta) + C_{\hat{x}}\delta + C_{\hat{x}}({\rm
e}^{-\lambda \delta}-1).
\end{array}
\end{equation*}

Therefore
\begin{equation*}
\begin{array}{r@{\ }c@{\ }l}
\| \Phi(\hat{x})-\Phi(S(\delta)\hat{x}) \| &\leq& \|
\Phi(\hat{x})-\Phi(y_{\hat{x}}(\delta)) \| +
L(\Phi)C_{\hat{x}}\delta \\[.7pc]
&\leq& \small{o}(\delta) + C \delta + C({\rm e}^{-\lambda
\delta}-1).
\end{array}
\end{equation*}

This proves the claim that $D_A^-\Phi(\hat{x}) <\infty$ and hence
from Lemma~\ref{c4.9}, it follows that as $t \to 0$,
\begin{equation} \label{3}
\left. \begin{array}{@{\hskip -4pc}r@{\ }c@{\ }l@{}} {\rm
e}^{-\lambda t} \Phi(y_{\hat{x}}(t)) - \Phi(\hat{x}) & \leq
&\displaystyle
\frac{{\rm e}^{-\lambda t}-1}{\lambda} D_A^- \Phi (\hat{x})\\[1.3pc]
&&\displaystyle + L(\psi) \int_0^t {\rm e}^{-\lambda s} \|
f(y_{\hat{x}}(s),u^{1,t}(s),d^1,{\beta}^t(u^{1,t}(\cdot),d^1)(s))
\|~{\rm d}s \\[1.3pc]
&&\displaystyle+\int_0^t {\rm e}^{-\lambda s} [\langle
D\phi(y_{\hat{x}}(s)),
f(y_{\hat{x}}(s),u^{1,t}(s),d^1,\\[1.5pc]
&&\quad\, {\beta}^t(u^{1,t}(\cdot),d^1)(s)) \rangle - \lambda
\Phi(y_{\hat{x}}(s))]~{\rm d}s + \small{o}(t).
\end{array}\!\right\}
\end{equation}

Combining (\ref{1}), (\ref{2}) and (\ref{3}), we obtain
\begin{align*}
\begin{array}{r@{\ }c@{\ }l}
\left[\dfrac{{\rm e}^{-\lambda t} \!-\!1}{\lambda}\right]
\dfrac{r}{2}
& \geq &{\rm e}^{-\lambda t} \Phi(y_{\hat{x}}(t))-\Phi(\hat{x}) \\
&&\displaystyle +\int_0^t {\rm e}^{-\lambda s}
k(y_{\hat{x}}(s),u^{1,t}(s),d^1,{\beta}^t(u^{1,t}(\cdot),d^1)(s))~{\rm
d}s \!+\!\small{o}(t) \\[1pc]
&\geq& \small{o}(t).
\end{array}
\end{align*}

This contradiction proves that $V_-$ is an approximate subsolution
of (HJI1$+$).

To prove that $V_-$ is an approximate supersolution of (HJI1$+$),
let $\hat{x} \in B_R(0)$ be a local minimum of $V^{d^1,d^2}_-
-\Phi$. Without any loss of generality, we may assume that
$V^{d^1,d^2}_- (\hat{x}) = \Phi(\hat{x})$. If $V^{d^1,d^2}_-
(\hat{x}) =M^{d^1,d^2}_-[V_-](\hat{x})$ or $V^{d^1,d^2}_-
(\hat{x}) =N[V_-^{d^1,d^2}](\hat{x})$, then we are done. Assume
that $V^{d^1,d^2}_- (\hat{x}) < \min
(M^{d^1,d^2}_-[V_-](\hat{x}),N[V_-^{d^1,d^2}](\hat{x}))$. In this
case, we need to show that
\begin{equation*}
\lambda \Phi(x) + D_A^+\Phi(\hat{x}) +H^{d^1,d^2}_+(\hat{x},
D\phi(\hat{x})) +C_R L(\psi) =: \hat{r} \geq 0.
\end{equation*}

If possible, let $\hat{r} <0$. Then
\begin{align*}
&\lambda \Phi(\hat{x}) + D_A^+\Phi(\hat{x})- \langle
D\phi(\hat{x}),f(\hat{x},\bar{u}^1,d^1,u^2,d^2) \rangle\\[.4pc]
&\quad\, - k((\hat{x},\bar{u}^1,d^1,u^2,d^2) + C_R L(\psi) \leq
\frac{\hat{r}}{2},
\end{align*}
for some $\bar{u}^1 \in U^1$ and all $u^2 \in U^2$. This implies
that, for all $\beta \in \Delta^{d^2}_0$ and $s$ small enough
\begin{align*}
&\lambda \Phi(y_{\hat{x}}(s)) + D_A^+\Phi(\hat{x}) - \langle
D\phi(y_{}(s)),
f(y_{\hat{x}}(s),\bar{u}^1,d^1,\beta(\bar{u}^1,d^1)(s)) \rangle\\[.3pc]
&\quad\, - k(y_{\hat{x}}(s),\bar{u}^1,d^1,\beta(\bar{u}^1,d^1)(s))
+ C_R L(\psi) \leq \frac{\hat{r}}{4}.
\end{align*}

Multiplying throughout by ${\rm e}^{-\lambda s}$ and integrating
from $0$ to $t$, we get
\begin{equation} \label{4}
\left.\begin{array}{@{\hskip -4.7pc}ll} &\left\{
D_A^+\Phi(\hat{x}) + C_R L(\psi) -\dfrac{\hat{r}}{4} \right\} \,
\left[\dfrac{{\rm
e}^{-\lambda t}-1}{-\lambda}\right]\\[1pc]
&\quad\, \displaystyle + \int_0^t {\rm e}^{-\lambda s}[\lambda
\Phi(y_{\hat{x}}(s))- \langle D\phi(y_{}(s)),
f(y_{\hat{x}}(s),\bar{u}^1,d^1,\beta(\bar{u}^1,d^1)(s)) \rangle ]~{\rm d}s\\[1.5pc]
&\quad\, \displaystyle - \int_0^t {\rm e}^{-\lambda s}
k(y_{\hat{x}}(s),\bar{u}^1,d^1,\beta(\bar{u}^1,d^1)(s))~{\rm d}s
\, \leq \, 0
\end{array}\right\}.
\end{equation}

Now by Lemma~\ref{s2i3}(iv), for $t$ small enough, there exists
$\beta^t \in \Delta^{d^2}_0$ such that
\begin{align} \label{5}
&V^{d^1,d^2}_- (\hat{x}) + \small{o}(t)\nonumber\\[.4pc]
&\quad\, \geq \int_0^t {\rm e}^{-\lambda s}
k(y_{\hat{x}}(s),\bar{u}^1,d^1,\beta^t(\bar{u}^1,d^1)(s))~{\rm d}s
+ {\rm e}^{-\lambda t} V^{d^1,d^2}_- (y_{\hat{x}}(t)) .
\end{align}

Now proceeding as in the first part of the theorem, we can show
that $D_A^+ \Phi(\hat{x}) > - \infty$. Using Lemma \ref{c4.9} for
$-\Phi$, we get
\begin{align}\label{6}
\left. \begin{array} {@{\hskip -4pc}r@{\ }c@{\ }l@{}} -{\rm
e}^{-\lambda t} \Phi(y_{\hat{x}}(t)) + \Phi(\hat{x}) & \leq &
\displaystyle\frac{{\rm e}^{-\lambda t}-1}{-\lambda} D_A^+ \Phi (\hat{x})\\[1.5pc]
&+&\displaystyle L(\psi) \int_0^t {\rm e}^{-\lambda s} \|
f(y_{\hat{x}}(s),u^{1,t}(s),d^1,{\beta}^t(u^{1,t}(\cdot),d^1)(s))
\|~{\rm d}s\\[1.5pc]
&-&\displaystyle \int_0^t {\rm e}^{-\lambda s} [\langle
D\phi(y_{\hat{x}}(s)),
f(y_{\hat{x}}(s),u^{1,t}(s),d^1,\\[1.5pc]
&&{\beta}^t(u^{1,t}(\cdot),d^1)(s)) \rangle -\displaystyle \lambda
\Phi(y_{\hat{x}}(s))]~{\rm d}s + \small{o}(t).
\end{array}\right\}\hskip -1pc\phantom{0}
\end{align}

From (\ref{4}), (\ref{5}) and (\ref{6}), we get
\begin{equation*}
\frac{\hat{r}}{4} \left[\frac{{\rm e}^{-\lambda
t}-1}{\lambda}\right] \leq \small{o}(t).
\end{equation*}

This is a contradiction and proves the fact that $V_-$ is an
approximate supersolution of (HJI1$+$).\hfill $\blacksquare$
\end{proof}

In a similar fashion, we can show that $V_-$ is an approximate
viscosity solution of (HJI2$+$). Hence $V_-$ is an approximate
viscosity solution of (HJI$+$).

\section{Existence of value}

\setcounter{equation}{0} \setcounter{theore}{0}

\looseness 1 In this section we first prove the uniqueness of
solutions of (HJI$+$) and (HJI$-$). Next we prove that uniqueness
result for (HJI$+$) and (HJI$-$) holds true even if one is a
viscosity solution and other is an approximate viscosity solution.
Finally under Isaacs' type min--max condition we show that game
has a value proving the main theorem of the paper. Now we state
and prove two lemmas needed in the proof of uniqueness of (HJI$+$)
and (HJI$-$).

\begin{lem}\label{s3i1}
Assume {\rm (A2)}. Let $w$ be uniformly continuous. If
\begin{equation*}
w^{d^1 , d^2 }( y_0 ) = N [w ^{d^1 , d^2 }] (y_0) = w^{d^1 , d^2
}( y_0 + \xi_0 ) + l(\xi_0),
\end{equation*}
then there exists $\sigma > 0$ {\rm (}which depends only on $w)$
such that for all $ y \in \bar {B} ( y_0 + \xi_0, \; \sigma ), $
\begin{equation*}
w^{d^1 , d^2 }( y ) < N [w ^{d^1 , d^2 }] (y).
\end{equation*}
\end{lem}

\begin{proof}
The proof closely mimics the corresponding result in the finite
dimensional case \cite{Y3}. We however prove it for the sake of
completeness. Let
\begin{equation*}
w^{d^1 , d^2 }( y_0 ) = N [w ^{d^1 , d^2 }] (y_0) = w^{d^1 , d^2
}( y_0 + \xi_0 ) + l(\xi_0).
\end{equation*}

Then, for every $ \xi_1 \in K $,
\begin{align*}
\hskip -4pc w^{d^1 , d^2 }( y_0 + \xi_0 + \xi_1 ) + l(\xi_1 ) -
w^{d^1 , d^2 }( y_0 + \xi_0) &= w^{d^1 , d^2 }( y_0 + \xi_0 +
\xi_1 ) + l(\xi_1)\\[.3pc]
&\quad\, - w^{d^1 , d^2 }( y_0 ) + l(\xi_0)\\[.3pc]
&\geq - l(\xi_0 + \xi_1 ) + l(\xi_0) + l(\xi_1).
\end{align*}
Now,
\begin{equation*}
N [w ^{d^1 , d^2 }] (y_0 + \xi_0 ) = \inf_{ \xi_1 \in K } w^{d^1 ,
d^2 }( y_0 + \xi_0 + \xi_1 ) + l(\xi_1).
\end{equation*}
By (\ref{s1e3}), this infimum will be attained in some $R$ ball.
Hence we can write
\begin{equation*}
N [w ^{d^1 , d^2 }] (y_0 + \xi_0 ) = \inf_{ \xi_1 \in K, \,
|\xi_1| < R } w^{d^1 , d^2 }( y_0 + \xi_0 + \xi_1 ) + l(\xi_1).
\end{equation*}
Now, taking infimum over $ |\xi_1| < R$ on both sides of the
earlier inequality we will\break have
\begin{align*}
\hskip -4pc N [ w^{d^1 , d^2 }( y_0 + \xi_0 )] - w^{d^1 , d^2 }(
y_0 + \xi_0 ) &\geq \inf\limits_{\xi_1 \in K, \, |\xi_1|<R } [
l(\xi_0) + l(\xi_1)- l(\xi_0 + \xi_1 )] \\[.3pc]
&= \bar l > 0.
\end{align*}
By using uniform continuity of $ w^{d^1 , d^2 }$ and $ N [ w^{d^1,
d^2 }]$ we get a $ \sigma $ such that for all $ y \in \bar {B} (
y_0 + \xi_0, \; \sigma ),$
\begin{equation*}
w^{d^1 , d^2 }( y ) < N [w ^{d^1 , d^2 }] (y).
\end{equation*}

$\left.\right.$\vspace{-2.8pc}

\hfill $\blacksquare$
\end{proof}

\begin{lem}\label{s3i2}
Assume {\rm (A1)} and {\rm (A2)}. Then the following results are
true.\vspace{-.3pc}

\begin{enumerate}
\renewcommand\labelenumi{\rm (\roman{enumi})}
\leftskip .2pc
\item Any supersolution $w$ of {\rm (HJI1}$+${\rm )} satisfies $w^{d^{1},d^{2}}
\geq M^{d^1,d^2}_+[w]$ for all $d^1,d^2$.

\item Any subsolution $w$ of {\rm (HJI2}$+${\rm )} satisfies $w^{d^{1},d^{2}}
\leq \min(M^{d^1,d^2}_-[w], N[w^{d^1,d^2}])$ for all
$d^1,d^2$.\vspace{-.3pc}
\end{enumerate}
\end{lem}

\begin{proof}
Let $w$ be a supersolution of (HJI1$+$). If possible, let
\begin{equation*}
w^{d^1,d^2}(x_0) < M_+^{d^1,d^2}[w](x_0).
\end{equation*}

By continuity, the above holds for all $x$ in an open ball $B$
around $x_0$. As in Lemma~1.8(d), p.~30 in \cite{BD}, we can show
that there exists $y_0 \in B$ and a smooth map $\phi$ such that
$w^{d^1,d^2} - \phi$ has local minimum at $y_0$. Since $w$ is a
supersolution of (HJI1$+$), this will lead to
\begin{equation*}
w^{d^1,d^2}(y_0) \geq M_+^{d^1,d^2}[w](y_0),
\end{equation*}
a contradiction. This proves (i). The proof of (ii) is
similar.\hfill $\blacksquare$
\end{proof}

Next we present the proof of the uniqueness theorem.

\begin{theor}[\!] \label{s3i3}
Assume {\rm (A1)} and {\rm (A2)}. Let $v$ and $w \in {\rm
BUC}(\mathbb{E}; \mathbb{R}^{ m_1 \times m_2})$ be viscosity
solutions of {\rm (HJI}$+${\rm )} {\rm (}or {\rm (HJI}$-${\rm ))}.
Then $ v = w $.
\end{theor}

\begin{proof} We prove the uniqueness for (HJI$+$).
The result for (HJI$-$) is similar.

Let $v$ and $w$ be viscosity solutions of (HJI$+$). We prove
$v^{d^1, d^2} \leq w^{d^1, d^2} $ for all $ d^1, d^2$. In a
similar fashion we can prove that $ w^ {d^1, d^2} \leq v ^ {d^1,
d^2}$ for all $ d^1, d^2$.

For $ ( d^1, d^2 ) \in D^1 \times D^2 $, define $ \Phi^ {{d^1,
d^2}}{\rm :}\ \mathbb{E} \times \mathbb{E} \to \mathbb{R}$ by
\begin{equation*}
\Phi^ {{d^1, d^2}} ( x, y) = v^ {d^1, d^2} (x) - w ^ {d^1, d^2}
(y) - \frac{| x- y |^2} {2 \epsilon} - \kappa [\langle
x\rangle^{\bar m} + \langle y\rangle^{\bar m}],
\end{equation*}
where $\bar m \in ( 0, 1) \cap \Big( 0, \frac{ \lambda}{\| f
\|_{\infty} }\Big)$ is fixed, $\kappa ,\epsilon \in (0, 1) $ are
parameters, and $\langle x \rangle^ {\bar m} = ( 1 + \|x\|^2 )^
{{\bar m }/2} $.

Note that $x \rightarrow \langle x \rangle^{\bar m} $ is Lipschitz
continuous with Lipschitz constant 1.

We first fix $\kappa > 0$. Let $ ( x_\epsilon, y_\epsilon)$ and $
( d^1_\epsilon, d^2_\epsilon ) $ be such that
\begin{equation*}
\Phi^{d^1_\epsilon , d^2_\epsilon } ( x_\epsilon, y_\epsilon) =
\sup \limits_{x,y} \max \limits_{d^1, d^2} \Phi^{d^1, d^2} ( x, y)
- \epsilon [d(x,x_{\epsilon}) + d(y,y_{\epsilon})]
\end{equation*}
and
\begin{equation*}
\Phi^{d^1_\epsilon , d^2_\epsilon } ( x_\epsilon, y_\epsilon) \geq
\sup \limits_{x,y} \max \limits_{d^1, d^2} \Phi^{d^1, d^2} ( x, y)
-\epsilon.
\end{equation*}

Here $d(\cdot , \cdot)$ is the Tataru's distance defined by
\begin{equation*}
d(x,y)=\inf_{t\geq 0}~[ t+ \| x-S(t)y \|]
\end{equation*}
(see \cite{CLVI} for more details). By Lemma~\ref{s3i2}, we have
\begin{align}\label{s3e1}
w^{d_\epsilon^1 , d_\epsilon^2 }( y_\epsilon ) \leq N [w
^{d^1_\epsilon , d^2_\epsilon }] (y_\epsilon),\\[.3pc]\label{s3e2}
w^{d_\epsilon^1 , d_\epsilon^2 }( y_\epsilon )\leq
{M_-}^{d^1_\epsilon , d^2_\epsilon } [w]
(y_\epsilon),\\[.3pc]\label{s3e3} v^{d_\epsilon^1 , d_\epsilon^2 }(
x_\epsilon )\geq {M_+}^ {d^1_\epsilon , d^2_\epsilon } [v]
(x_\epsilon) .
\end{align}

If we have strict inequality in all the above three inequalities
(that is, (\ref{s3e1}), (\ref{s3e2}) and (\ref{s3e3})), then by
the definition of viscosity sub and super solutions we will have
\begin{align}\label{h1}
\lambda v^{ d^1_\epsilon, d^2_\epsilon} ( x_\epsilon ) + D_A^-
\Phi_1(x_{\epsilon}) + H_{+,\underline{\epsilon}} ^{ d^1_\epsilon,
d^2_\epsilon} \; \left( x_\epsilon, \frac{x_\epsilon -
y_\epsilon}{ \epsilon} + \kappa \bar m  \langle x_\epsilon
\rangle^{ \bar {m} -2} x_\epsilon \right) \leq 0, \\ \label{h2}
\lambda w^{ d^1_\epsilon, d^2_\epsilon} ( y_\epsilon ) +
 D_A^+ \Phi_2(y_{\epsilon}) +
H_{+,\bar{\epsilon}}^{ d^1_\epsilon, d^2_\epsilon} \left(
y_\epsilon, \frac{x_\epsilon - y_\epsilon}{ \epsilon} - \kappa
\bar m \langle y_\epsilon \rangle^{ \bar {m} -2} y_\epsilon\right)
\geq 0,
\end{align}
where
\begin{equation*}
\begin{array} {rcl}
\Phi_1(x) &=& \dfrac{| x-y_{\epsilon}|^2}{2\epsilon} +
\kappa \langle x \rangle^{\bar{m}} + \, \epsilon d(x,x_{\epsilon}),\\[.9pc]
-\Phi_2(y) &=& \dfrac{| x_{\epsilon} -y|^2}{2\epsilon} + \kappa
\langle y \rangle^{\bar{m}} + \, \epsilon d(y,y_{\epsilon}).
\end{array}
\end{equation*}

In this case we can proceed by the usual comparison principle
method as in \cite{CLVI}.

Therefore it is enough to show that for a proper auxiliary
function strict inequality occurs in (\ref{s3e1}), (\ref{s3e2})
and (\ref{s3e3}) at the maximizer. We achieve this in the
following three steps.

\begin{step}{\rm
There are two cases to consider; either there is no sequence
$\epsilon_n \downarrow 0 $ such that strict inequality holds in
(\ref{s3e1}) or there is some sequence for which strict inequality
holds in (\ref{s3e2}).

If there is no sequence $\epsilon_n \downarrow 0$ such that
\begin{equation} \label{s3e4a}
w^{d_{\epsilon_n}^1 , d_{\epsilon_n}^2 }( y_{\epsilon_n} ) < N [w
^{d^1_{\epsilon_n} , d^2_{\epsilon_n} }] (y_{\epsilon_n}) \mbox{
for all $n$},
\end{equation}
then we have equality in (\ref{s3e1}) for all $\epsilon$ in some
interval $(0,\epsilon_0)$. By the definition of $N$ and the
assumptions (A2), for each $\epsilon \in (0,\epsilon_0)$, there
exists $\xi_\epsilon \in K$ such that
\begin{equation*}
M_0 := \sup_{0< \epsilon < \epsilon_0}\| \xi_{\epsilon} \|
< \infty
\end{equation*}
and
\begin{equation*}
w^{d_\epsilon^1 , d_\epsilon^2 }( y_\epsilon ) = N [w
^{d^1_\epsilon , d^2_\epsilon }] (y_\epsilon) = w^{d_\epsilon^1 ,
d_\epsilon^2 }( y_\epsilon + \xi_\epsilon ) + l(\xi_\epsilon).
\end{equation*}

Then,
\begin{align*}
\Phi^{d^1_\epsilon , d^2_\epsilon } ( x_\epsilon + \xi_\epsilon ,
y_\epsilon + \xi_\epsilon ) &= v^{d_\epsilon^1 , d_\epsilon^2 }(
x_\epsilon + \xi_\epsilon )- w^{d_\epsilon^1 , d_\epsilon^2 }(
y_\epsilon + \xi_\epsilon )- \frac{| x_\epsilon - y_\epsilon |^2} {2
\epsilon}\\[.5pc]
&\quad\, - \kappa [ \langle x_\epsilon + \xi_\epsilon \rangle^
{\bar m} +
\langle y_\epsilon + \xi_\epsilon \rangle^ {\bar m} ]
\end{align*}
\begin{align*}
&\geq v^{d_\epsilon^1 , d_\epsilon^2 }( x_\epsilon ) -
w^{d_\epsilon^1 , d_\epsilon^2 }( y_\epsilon )
- \frac{| x_\epsilon - y_\epsilon |^2} {2 \epsilon}\\[.5pc]
&\quad\, -\kappa [\langle x_\epsilon \rangle^{\bar{m}}+ \langle
y_\epsilon \rangle^{\bar{m}}] - 2 \kappa | \xi_\epsilon|\\[.5pc]
&= \Phi^{d^1_\epsilon , d^2_\epsilon } ( x_\epsilon , y_\epsilon)
- 2 \kappa | \xi_\epsilon | \\[.5pc]
&\geq \Phi^{d^1_\epsilon , d^2_\epsilon } ( x_\epsilon ,
y_\epsilon) - 2 \kappa M_0.
\end{align*}

Hence we have
\begin{equation}\label{s3e4}
\Phi^{d^1_\epsilon , d^2_\epsilon } ( x_\epsilon , y_\epsilon) -
\Phi^{d^1_\epsilon , d^2_\epsilon } ( x_\epsilon + \xi_\epsilon ,
y_\epsilon + \xi_\epsilon ) \leq 2 \kappa M_0.
\end{equation}

We will be using this difference to define the new auxiliary
function. Observe that $ |x_\epsilon| , |y_\epsilon| ,
|\xi_\epsilon| $ are bounded.

We define the new auxiliary function $\Psi$ by
\begin{align*}
\Psi^{d^1 , d^2 } (x, y) &= \Phi^{d^1 , d^2 } (x,y) -\epsilon
[d(x,x_\epsilon)+d(y,y_\epsilon)]\\[.5pc]
&\quad\, +2( 2 \kappa +\epsilon) M_0 \eta \left( \frac{x -
x_{\epsilon}-\xi_\epsilon}{\sigma}, \frac{y -
y_{\epsilon}-\xi_\epsilon}{\sigma}\right),
\end{align*}
where $\sigma = \sigma(w)$ is the constant coming from Lemma~3.1
and $\eta{\rm :}\ \mathbb{R}^d \times \mathbb{R}^d \rightarrow
\mathbb{R} $ is a smooth function with the following properties:
\begin{enumerate}
\renewcommand\labelenumi{\arabic{enumi}.}
\leftskip -.1pc
\item $\hbox{supp} ( \eta ) \subset B ((0,0), \; 1)$,
\item $ 0 \leq \eta \leq 1 $,
\item $\eta (0,0) =1 $ and $\eta < 1 $ if $(x,y) \neq (0,0)$,
\item $ | D \eta | \leq 1 $.
\end{enumerate}

Now by the definition of $\Psi$,
\begin{align*}
\Psi^{d_\epsilon^1 , d_\epsilon^2 } ( x_{\epsilon} + \xi_\epsilon,
y_{\epsilon} + \xi_\epsilon ) &= \Phi^{d_\epsilon^1 , d_\epsilon^2
}( x_{\epsilon} + \xi_\epsilon, y_{\epsilon} + \xi_\epsilon )\\[.4pc]
&\quad\, -\epsilon[d(x_\epsilon+
\xi_\epsilon,x_\epsilon)+d(y_\epsilon+\xi_\epsilon,
y_\epsilon)]\\[.4pc]
&\quad\, + 2(2 \kappa+ \epsilon) M_0\\[.4pc]
&\geq \Phi^{d_\epsilon^1 , d_\epsilon^2 } ( x_{\epsilon} ,
y_{\epsilon}) +2\kappa M_0
\end{align*}
and
\begin{align*}
\Psi ^{d_\epsilon^1 , d_\epsilon^2 } (x, y) &= \Phi^{d_\epsilon^1
, d_\epsilon^2 } (x, y) \;\; \quad \mbox{if}\;\; | x- x_{\epsilon}
|^2 + |y - y_{\epsilon} |^2 \geq \sigma^2.
\end{align*}

Hence $\Psi^{d_\epsilon^1 , d_\epsilon^2 } $ attains its maximum
in the $\sigma$ ball around $ ( x_{\epsilon} + \xi_\epsilon,
y_{\epsilon} + \xi_\epsilon )$ at (say) the point $
(\hat{x}_\epsilon, \hat{y}_\epsilon) .$ By Lemma (\ref{s3i1}), we
now know that, for all $ \epsilon \in (0, \epsilon_0)$,
\begin{equation*}
w^{d_\epsilon^1 , d_\epsilon^2 }( \hat {y}_\epsilon ) <
N [w ^{d^1_\epsilon , d^2_\epsilon }] (\hat {y}_\epsilon ).
\end{equation*}

If there is a sequence $\epsilon_n \downarrow 0$ such that
(\ref{s3e4a}) holds for all $n$, then, along this sequence, we
proceed to Step 2 with $\Psi^{d^1,d^2} = \Phi^{d^1,d^2}$ and
$(\hat{x}_{\epsilon_n},\hat{y}_{\epsilon_n})= ( {x}_{\epsilon_n},
{y}_{\epsilon_n})$.

Thus we always have a sequence $\epsilon_n \downarrow 0$ along
which we have
\begin{equation} \label{s3e5}
w^{d_{\epsilon_n}^1 , d_{\epsilon_n}^2 }( \hat{y}_{\epsilon_n} ) <
N [w ^{d^1_{\epsilon_n} , d^2_{\epsilon_n} }]
(\hat{y}_{\epsilon_n}),
\end{equation}
with $(\hat{x}_{\epsilon_n},\hat{y}_{\epsilon_n})$ being a
maximizer of $\Psi^{d^1_{\epsilon_n},d^2_{\epsilon_n}}$. Since
$D^1 \times D^2$ is a finite set, without any loss of generality,
we may assume that $(d_{\epsilon_n}^1 , d_{\epsilon_n}^2 )=
(d_{0}^1 , d_{0}^2 )$ for all $n$. In the next step we check what
happens to the inequality (\ref{s3e2}) at the maximizer of the new
auxiliary function $ \Psi^{d^1_{0},d^2_{0}}$.}
\end{step}

\begin{step}{\rm
Now for each fixed $n$, we proceed as follows. There are two cases
either $w^{d_{0}^1 , d_{0}^2 }( \hat{y}_{\epsilon_n} ) =
M_-^{d^1_{0} , d^2_{0} } [w] (\hat{y}_{\epsilon_n}) $ or
$w^{d_{0}^1 , d_{0}^2 }( \hat{y}_{\epsilon_n} ) < M_-^{d^1_{0} ,
d^2_{0} } [w] (\hat{y}_{\epsilon_n}) $.

If
\begin{equation*}
w^{d_{0}^1 , d_{0}^2 }( \hat{y}_{\epsilon_n} ) = M_-^{d^1_{0} ,
d^2_{0} } [w] (\hat{y}_{\epsilon_n}),
\end{equation*}
then by the definition of $M_-^{d^1_0 , d^2_0 } $, there exists $
d^2_{n_1} \in D^2$ such that
\begin{equation}\label{s3e6}
w^{d_0^1 , d_0^2 }( \hat{y}_{\epsilon_n} ) = w^{d_0^1 , d_{n_1}^2
}( \hat{y}_{\epsilon_n} ) + c^2 (d^2_0, d_{n_1}^2 ).
\end{equation}

We know that
\begin{equation*}
v^{d_0^1 , d_0^2 }( \hat{x}_{\epsilon_n} ) \leq v^{d_0^1 ,
d_{n_1}^2 }( \hat{x}_{\epsilon_n} ) + c^2 ( d^2_0, d_{n_1}^2 ).
\end{equation*}

Hence,
\begin{align*}
\Psi ^{d_0^1 , d_0^2 } ( \hat{x}_{\epsilon_n},
\hat{y}_{\epsilon_n} ) -\Psi ^{d_0^1 , d_{n_1}^2 } (
\hat{x}_{\epsilon_n}, \hat{y}_{\epsilon_n} ) &\leq v^{d_0^1 ,
d_{n_1}^2 }( \hat {x}_{\epsilon_n} ) + c^2 ( d^{2}_0, d_{n_1}^2)\\[.4pc]
&\quad\, - w^{d_0^1 , d_{n_1}^2 }( \hat {y}_{\epsilon_n} )
- c^2 ( d^2_0 , d_{n_1}^2 )\\
&\leq 0.
\end{align*}

Hence we get
\begin{equation*}
\Psi ^{d_0^1 , d_0^2 } ( \hat{x}_{\epsilon_n},
\hat{y}_{\epsilon_n} ) \leq \Psi ^{d_0^1 , d_{n_1}^2 } (
\hat{x}_{\epsilon_n}, \hat{y}_{\epsilon_n} ).
\end{equation*}

Now if strict inequality holds in $ w^{d_0^1 , d_{n_1}^2 } (
\hat{y}_{\epsilon_n} ) \leq M_-^{d_0^1 , d_{n_1}^2 } [w]
(\hat{y}_{\epsilon_n})$, then we are done; else we repeat the
above argument and get $ d_{n_2}^2 \in D^2 $ such that
\begin{equation*}
w^{d_0^1 , d_{n_1}^2 } ( \hat{y}_{\epsilon_n} )= w^{d_0^1 ,
d_{n_2}^2 } ( \hat{y}_{\epsilon_n} ) + c^2 ( d_{n_1}^2 , d_{n_2}^2
).
\end{equation*}

Now
\begin{align*}
w^{d_0^1 , d_{n_2}^2 } ( \hat{y}_{\epsilon_n} ) &= w^{d_0^1 ,
d_{n_1}^2 } ( \hat{y}_{\epsilon_n} )- c^2 ( d_{n_1}^2 , d_{n_2}^2 )\\[.3pc]
&\geq w^{d_0^1 , d_{n_1}^2 } ( \hat{y}_{\epsilon_n} )
- c^2_{0} \\[.3pc]
&\geq w^{d_0^1 , d_0^2 }( \hat{y}_{\epsilon_n} ) - 2 c^2_{0}.
\end{align*}

Proceeding in similar fashion, after finitely many steps,
boundedness of $w$ will be contradicted and hence for some $d^2_n
= d^2_{n_j} \in D^2 $, we must have
\begin{equation} \label{s3e7}
w^{d_0^1 , d_{n}^2 } ( \hat{y}_{\epsilon_n} ) < M_-^{d_0^1 ,
d_{n}^2 } [w] (\hat{y}_{\epsilon_n})
\end{equation}
and
\begin{equation}\label{*}
\Psi^{d^1_0,d^2_0}( \hat{x}_{\epsilon_n}, \hat{y}_{\epsilon_n})
\leq \Psi^{d^1_0,d^2_n}( \hat{x}_{\epsilon_n},
\hat{y}_{\epsilon_n}).
\end{equation}

On the other hand, if
\begin{equation*}
w^{d_{0}^1 , d_{0}^2 }( \hat{y}_{\epsilon_n} ) < M_-^{d^1_{0} ,
d^2_{0} } [w] (\hat{y}_{\epsilon_n}),
\end{equation*}
then we proceed by taking $d^2_n = d^2_0$.}
\end{step}

\begin{step}{\rm
For each fixed $n$, we proceed as follows: If
\begin{equation*}
v^{d_{0}^1 , d_{n}^2 }( \hat{x}_{\epsilon_n} ) = M_+^{d^1_{0} ,
d^2_{n} } [v] (\hat{x}_{\epsilon_n}),
\end{equation*}
then we proceed as in Step 2 and obtain $d^1_n=d^1_{n_i} \in D^1$
such that
\begin{equation*}
v^{d_{n}^1 , d_{n}^2 }( \hat{x}_{\epsilon_n} ) > M_+^{d^1_{n} ,
d^2_{n} } [v] (\hat{x}_{\epsilon_n})
\end{equation*}
and
\begin{equation}\label{**} \Psi^{d^1_0,d^2_n}(
\hat{x}_{\epsilon_n}, \hat{y}_{\epsilon_n}) \leq
\Psi^{d^1_n,d^2_n}( \hat{x}_{\epsilon_n}, \hat{y}_{\epsilon_n}).
\end{equation}

If
\begin{equation*}
v^{d_{0}^1 , d_{n}^2 }( \hat{x}_{\epsilon_n} ) = M_+^{d^1_{0} ,
d^2_{n} } [v] (\hat{x}_{\epsilon_n}),
\end{equation*}
then we proceed by taking $d^1_n=d^1_0$.

Thus, for every $n$, $(\hat{x}_{\epsilon_n},\hat{y}_{\epsilon_n})$
is a maximizer of $\Psi^{d^1_n,d^2_n}$ and
\begin{align}\label{S3e8}
w^{d_n^1 , d_n^2 }( \hat{y}_{\epsilon_n} ) < N [w ^{d^1_n , d^2_n
}] (\hat{y}_{\epsilon_n}), \\[.4pc]\label{s3e9}
 w^{d_n^1 , d_n^2 }(\hat{y}_{\epsilon_n} ) <
M_-^{d^1_n , d^2_n } [w] (\hat{y}_{\epsilon_n}),\\[.4pc]\label{s3e10}
 v^{d_n^1 , d_n^2 }( \hat{x}_{\epsilon_n} ) > M_+^
{d^1_n , d^2_n } [v] (\hat{x}_{\epsilon_n}).
\end{align}

Also by using (\ref{*}) and (\ref{**}), we have that
\begin{equation}\label{***}
\Psi^{d^1_0,d^2_0}( \hat{x}_{\epsilon_n}, \hat{y}_{\epsilon_n})
\leq \Psi^{d^1_n,d^2_n}( \hat{x}_{\epsilon_n},
\hat{y}_{\epsilon_n}).
\end{equation}

Now we define test functions $ \phi_{1n}$ and $\phi_{2n}$ as
follows:
\begin{align*}
\phi_{1n} (x) &= w^{d_{n}^1 , d_{n}^2 }( \hat {y}_{\epsilon_n} ) +
\frac{| x- \hat{y}_{\epsilon_n} |^2} {2 \epsilon_n} + \kappa [
\langle x \rangle  ^ {\bar m} + \langle
\hat{y}_{\epsilon_n}\rangle ^ {\bar m} ]\\[.3pc]
&\quad\, -2( 2 \kappa +\epsilon_n) M_0 \eta \left( \frac{ x-
x_{\epsilon_n}-\xi_{\epsilon_n}}{\sigma}, \frac{
\hat{y}_{\epsilon_n} -
y_{\epsilon_n}-\xi_{\epsilon_n}}{\sigma} \right), \\[.3pc]
\phi_{2n} ( y) &= v^{d_{n}^1 , d_{n}^2 }( \hat {x}_{\epsilon_n} )
 - \frac{| \hat{x}_{\epsilon_n}- {y} |^2} {2 \epsilon_n} - \kappa
[\langle \hat{x}_{\epsilon_n}\rangle  ^ {\bar m} + \langle
{y}\rangle  ^ {\bar m} ]\\[.3pc]
&\quad\, + 2 (2\kappa +\epsilon_n ) M_0  \eta \left(
\frac{\hat{x}_{\epsilon_n} - x_{\epsilon_n} - \xi_{\epsilon_n}
}{\sigma}, \frac{ {y} - y_{\epsilon_n} - \xi_{\epsilon_n}
}{\sigma}\right).
\end{align*}

Observe that
\begin{align*}
D \phi_{1n} ( \hat{x}_{\epsilon_n} ) &= \frac{\hat{x}_{\epsilon_n}
- \hat{y}_{\epsilon_n}}{\epsilon_n} + \kappa \bar{m} \langle
\hat{x}_{\epsilon_n} \rangle ^ {m-2} \hat{x}_{\epsilon_n}\\[.4pc]
&\quad\, - \frac{2(2 \kappa +\epsilon_n) M_0} {\sigma} \, D_x \eta
\left( \frac{\hat{x}_{\epsilon_n} - x_{\epsilon_n}-
\xi_{\epsilon_n}} {\sigma}, \frac{ \hat{y}_{\epsilon_n} - y_{\epsilon_n}-\xi_{\epsilon_n}}{\sigma}\right),\\
D \phi_{2n} ( \hat{y}_{\epsilon_n} ) &= \frac{\hat{x}_{\epsilon_n}
- \hat{y}_{\epsilon_n}}{\epsilon_n} - \kappa \bar{m} \langle
\hat{y}_{\epsilon_n} \rangle ^ {m-2} \hat{y }_{\epsilon_n}\\[.4pc]
&\quad\, + \frac{ 2 (2\kappa +\epsilon_n) M_0}{\sigma} \, D_y \eta
\left( \frac{\hat{x}_{\epsilon_n} -
x_{\epsilon_n}-\xi_{\epsilon_n}}{\sigma}, \frac{
\hat{y}_{\epsilon_n} - y_{\epsilon_n}-\xi_{\epsilon_n}}
{\sigma}\right).
\end{align*}

Note that $ v^{d_{n}^1 , d_{n}^2 }- \phi_{1n}$ attains its maximum
at $ \hat{x}_{\epsilon_n}$ and $ w^{d_{n}^1 , d_{n}^2 } -
\phi_{2n}$ attains its minimum at $ \hat {y}_{\epsilon_n} $.
Hence, as $\epsilon_n \to 0$, we have
\begin{align*}
\hskip -4pc\lambda [ v^{ d^1_{n}, d^2_n} ( \hat {x}_{\epsilon_n} )
\!-\! w^{ d^1_{n}, d^2_{n}} ( \hat {y}_{\epsilon_n} ) ] &\leq
H_+^{ d^1_{n}, d^2_{n} } ( \hat{y}_{\epsilon_n}, D \phi_{2n} (
\hat{y}_{\epsilon_n}))\!-\! H_+^{ d^1_{n}, d^2_{n}} (
\hat{x}_{\epsilon_n}, D \phi_{1n} ( \hat{x}_{\epsilon_n}) ) \!+\!
\small{o}(1)\\[.5pc]
&\leq L | \hat{x}_{\epsilon_n} - \hat{y}_{\epsilon_n}| \left( 1
+\left|\frac{\hat{x}_{\epsilon_n} -
\hat{y}_{\epsilon_n}}{\epsilon_n}\right|\right)
+ \small{o}(1) \\[.5pc]
&\quad\, + \| f \|_\infty \left[\kappa \bar{m} (\langle
\hat{x}_{\epsilon_n} \rangle ^{\bar{m}-1} + \langle
\hat{y}_{\epsilon_n} \rangle ^{\bar{m}-1}) +\frac{ 4
(2\kappa+\epsilon) M_0} {\sigma} \right].
\end{align*}

Note that we have used $ | {y}| \leq ( 1 + | {y}|^2 ) ^ {1/2} =
\langle {y} \rangle$ to get the above inequality. Now as $ \bar
{m} -1 < 0 $, it follows that
\begin{align*}
v^{ d^1_{n}, d^2_{n}} ( \hat{x}_{\epsilon_n}) - w^{ d^1_{n},
d^2_{n}} ( \hat{y}_{\epsilon_n} ) &\leq \frac {L} {\lambda} |
\hat{x}_{\epsilon_n} - \hat{y}_{\epsilon_n}| + \frac{L} {\lambda}
\frac {{|\hat{x}_{\epsilon_n}- \hat{y}_{\epsilon_n}|}^2}
{\epsilon_n} + \small{o}(1)\\[.5pc]
&\quad\, + \frac { 2 \| f \|_\infty \kappa \bar{m}} {\lambda} +
\frac {4 \| f \|_\infty \kappa M_0} {\sigma \lambda}\\[.5pc]
&\leq \frac { 2 \| f \|_\infty \kappa \bar{m}} {\lambda} +
 \frac {4 \| f \|_\infty \kappa M_0} {\sigma \lambda}
+ o(1) ; \;\mbox{as}\; \epsilon_n \downarrow 0.
\end{align*}

For any $ x \in \mathbb{R}^d$ and $( d^1, d^2 ) \in D^1 \times
D^2, $
\begin{align*}
\hskip -4pc v^{ d^1, d^2} (x) - w^{d^1, d^2} (x) - 2 \kappa
\langle x \rangle^{\bar{m}} &= \Phi^{d^1, d^2} (x,x) \\[.5pc]
&\leq \Phi^{d^1_0, d^2_0} (x_{\epsilon_n}, y_{\epsilon_n} )
+ \epsilon_n [d(x,x_{\epsilon_n})+d(x,y_{\epsilon_n}) ]\\[.5pc]
&\leq \Psi ^{d^1_0, d^2_0 }( x_{\epsilon_n}, y_{\epsilon_n} ) +
o(1) ; \;\mbox{as}\; \epsilon_n \downarrow 0.
\end{align*}

Since $\Psi^{d^1_0, d^2_0}$ attains its maximum at
$(\hat{x}_{\epsilon_n}, \hat{y}_{\epsilon_n} ) $ in the $\sigma$
ball around $ (x_{\epsilon_n} + \xi_{\epsilon_n}, y_{\epsilon_n} +
\xi_{\epsilon_n} )$, for an appropriate constant $C$, we have
\begin{align*}
\Psi ^{d^1_0, d^2_0 }( x_{\epsilon_n}, y_{\epsilon_n} ) &\leq
\Psi^{d^1_0, d^2_0 }( \hat{x}_{\epsilon_n},
\hat{y}_{\epsilon_n})\\[.5pc]
&\leq \Psi^{ d^1_{n}, d^2_{n} }
( \hat{x}_{\epsilon_n}, \hat{y}_{\epsilon_n}) \;\; \text{using}\; (\ref{***})
\end{align*}
\begin{align*}
&\leq v^{d^1_{n}, d^2_{n}} ( \hat{x}_{\epsilon_n}) - w^{d^1_{n},
d^2_{n}} ( \hat{y}_{\epsilon_n}) + C \kappa + o(1); \;\mbox{as}\;
\epsilon_n \downarrow 0 \\[.5pc]
&\leq \frac { 2 \| f \|_\infty \kappa \bar{m}} {\lambda} + \frac
{4 \| f \|_\infty \kappa M_0} {\sigma \lambda} + C \kappa + o(1);
\;\mbox{as}\; \epsilon_n \downarrow 0.
\end{align*}

Hence we get
\begin{align*}
&v^{ d^1, d^2} (x) - w^{d^1, d^2} (x) - 2 \kappa \langle x
\rangle^{\bar{m}}\\[.5pc]
&\quad\, \leq \frac { 2 \| f \|_\infty \kappa \bar{m}} {\lambda} +
\frac {4 \| f \|_\infty \kappa M_0} {\sigma \lambda} + C \kappa +
o(1); \;\mbox{as}\; \epsilon_n \downarrow 0.
\end{align*}
Now let $\epsilon_n \downarrow 0$ and then $\kappa \downarrow 0$,
to obtain
\begin{equation*}
v^{ d^1, d^2} (x) - w^{d^1, d^2} (x) \leq 0 .
\end{equation*}
This completes the proof of uniqueness for (HJI$+$).}\hfill
$\blacksquare$\vspace{-2pc}
\end{step}
\end{proof}

The above uniqueness result holds true if one is the viscosity
solution and the other is an approximate viscosity solution. This
is the content of the next theorem.

\begin{theor}[\!] \label{s3i4}
Assume {\rm (A1)} and {\rm (A2)}. Let $ v $ and $ w \in {\rm BUC}
(\mathbb{E}; \mathbb{R}^{ m_1 \times m_2}) $. Let $v$ be a
viscosity solution of {\rm (HJI}$+${\rm )} {\rm (}resp. {\rm
(HJI}$-${\rm ))} and $w$ be an approximate viscosity solution of
{\rm (HJI}$+${\rm )} {\rm (}resp. {\rm (}HJI$-${\rm ))}. Then $ v
= w $.
\end{theor}

\begin{proof}
The proof is similar to that of the previous theorem. The only
change here is in (\ref{h2}). Since $w$ is an approximate
viscosity solution, one gets
\begin{align*}
&\lambda w^{ d^1_\epsilon, d^2_\epsilon} ( y_\epsilon ) + D_A^+
\Phi_2(y_{\epsilon}) + H_{+}^{ d^1_\epsilon, d^2_\epsilon} \;
\left( y_\epsilon, \frac{x_\epsilon - y_\epsilon}{ \epsilon} -
\kappa \bar m \langle y_\epsilon \rangle ^{ \bar {m} -2}
y_\epsilon\right)\\[.5pc]
&\quad\, \geq - \epsilon \, C_R,
\end{align*}
where $R= \sup_{\epsilon >0} \| y_\epsilon \|$. Note that, for
fixed $\kappa$, $R < \infty$ and hence
\begin{equation*}
\epsilon C_R = \small{o}(1) \, \;\; \mbox{ as } \epsilon \to 0.
\end{equation*}

Once we have this inequality (instead of (\ref{h2})), we mimic all
other arguments in the proof of the previous theorem.\hfill
$\blacksquare$
\end{proof}

Now we can prove our main result stated in \S1, namely
Theorem~\ref{s1i1}.

\setcounter{section}{1}
\begin{pot}{\rm
Under the Isaacs min--max condition, (HJI$-$) and (HJI$+$)
coincide. Let us denote this equation by (HJI). As in \cite{CLVI},
by Perron's method, we can prove an existence of a viscosity
solution for (HJI) in $\hbox{BUC}({\mathbb{E}},{\mathbb{R}}^{m_1
\times m_2})$, the class of bounded uniformly continuous
functions. Let $W$ be any such viscosity solution. Now, by
Theorem~\ref{s2i6} we know that lower and upper value functions,
$V_- $ and $ V_+$ are approximate viscosity solutions of (HJI).
Therefore, by Theorem~\ref{s3i4}, $V_- = W = V_+$. This proves the
main result.}\hfill $\blacksquare$
\end{pot}

\setcounter{section}{3}

\section{Conclusions}

\setcounter{equation}{0} \setcounter{theore}{0}

We have studied two-person zero-sum differential games with hybrid
controls in infinite dimension. The minimizing player uses
continuous, switching, and impulse controls whereas the maximizing
player uses continuous and switching controls. The dynamic
programming principle for lower and upper value functions is
proved and using this we have established the existence and
uniqueness of the value under Isaacs min--max condition.

For finite dimensional problems, similar result has been obtained
by Yong \cite{Y3} under two additional assumptions:\vspace{.5pc}

\noindent (Y1) {\it Cheaper switching cost condition}
\begin{equation*}
\min_{\bar{d}^2 \neq d^2} c^2(d^2,\bar{d}^2)=:c^2_0 < l_0
=\inf_{\xi \in K}l(\xi).
\end{equation*}

\noindent (Y2) {\it Nonzero loop switching cost condition}

\hskip 1pc For any loop $\{ (d^1_i,d^2_i) \}_{i=1}^{j} \subset D^1
\times D^2 $, with the property that
\begin{align*}
&j \leq m_1 m_2, \, d^1_{j+1} =d^1_1, \, d^2_{j+1} =d^2_1 ;\\[1mm]
&\mbox{ either }\;\; d^1_{i+1} =d^1_i, \; \mbox{ or }
\;\;d^2_{i+1} =d^2_i \, \, \, \, \; \; \forall \,1 \leq i \leq j.
\end{align*}
It holds that
\begin{equation*}
\sum_{i=1}^{j} c^1(d^1_i,d^1_{i+1}) -
\sum_{i=1}^{j} c^2(d^2_i,d^2_{i+1}) \neq 0.
\end{equation*}

Thus our result not only extends the work of \cite{Y3} to infinite
dimensions but also proves the uniqueness of the viscosity
solutions of upper and lower SQVI without the above two conditions
(Y1) and (Y2). Also we have shown that under Isaacs' min--max
condition, the game has a value. Moreover, we have given explicit
formulation of dynamic programming principle for hybrid
differential games and have also proved it which is not done in
\cite{Y3}.

\section*{Acknowledgement}

The authors wish to thank M~K~Ghosh for suggesting the problem and
for several useful discussions. They also thank M~K~Ghosh and
Mythily Ramaswamy for carefully reading the manuscript and for
useful suggestions. Financial support from NBHM is gratefully
acknowledged.

\end{document}